%% file: lectures.tex
%

\input amstex
\documentstyle{ias-pcms} 
\input macros

\input boxedeps.tex
\magnification=1150
\voffset=-1cm

\SetepsfEPSFSpecial 
\HideDisplacementBoxes


\Volume{??}                	
\Year{1998}			 	

\SeriesTitle{Geometric Methods in Representation Theory}			
\ShortSeriesTitle{Geometric Methods}		
\Author{Kari Vilonen}			
\ShortAuthor{K.Vilonen}			
\AuthorAddress{Department of Mathematics, Brandeis University, Waltham,
MA 02254, USA}	
\AuthorEmail{vilonen\@math.brandeis.edu}			%



\EndTopInfo 

\head Introduction\endhead

The goal of this series of lectures is to survey and provide background for recent
joint work with Wilfried Schmid. This work  has appeared as a series of papers
\cite{SV1,SV2,SV3,SV4}. The type of geometric methods we discuss here were first
introduced to representation theory by Beilinson and Bernstein in \cite{BB1}, to solve
the Kazhdan-Lusztig conjectures. The localization technique of \cite{BB1} can be
used to translate questions in representation theory to questions about geometry of
complex algebraic varieties. Later, in \cite{K2}, Kashiwara initiated a research
program, as a series of conjectures, which extends the Beilinson-Bernstein picture.
This program was carried out in  \cite{KSd} and \cite{MUV} and is explained here in
lectures 4 and 5.  To explain the joint work with Schmid, we begin by
introducing the main technical tool, the characteristic cycle construction, in lecture
6. The primary objective of these lectures, the character formula and the proof of the
Barbasch-Vogan conjecture, are explained in lectures 7 and 8, respectively. In the
first lecture we present an overview of the lecture series and lectures 2 and 3
provide the necessary background material on sheaf theory and homological algebra.


\subhead Acknowledgments\endsubhead

I would like to thank Markus Hunziker for his help in the preparation of the lectures
at Park City, for drawing the pictures, and for writing the appendix that appears at the
end of these lectures. The notes taken at Park City by Upendra Kulkarni proved to be a
valuable help in preparing these notes. The final version of these notes were written
when the author was a guest at the Max Planck institute in Bonn. The author was
supported by grants from the NSF and NSA, and by a Guggenheim  fellowship.


\Lecture1{Overview}			

In this first lecture we give, in very rough terms, a basic outline of this
lecture series. Let $\GR$ be a semisimple
(linear, connected) Lie group. By a representation of $\GR$ we mean a
representation on a complete locally convex Hausdorff topological vector space, of
finite length. If the group $\GR$ acts on a manifold $X$ and $\cf$ is a sheaf on
$X$ which is $\GR$-equivariant, i.e., the group $\GR$ acts on the sections of the sheaf
$\cf$, then the cohomology groups $\oh^k(X,\cf)$ have a linear $\GR$-action. Let
$X$ denote the flag manifold of the complexification $G$ of $\GR$. Then the above
construction gives a functor:
$$
\{\text{certain $\GR$-equivariant sheaves on $X$}\} @>{\ \ \ \ \ }>>\{G_\Bbb R
\text{-representations}\}\,.
\tag1.1
$$
There is an analogous picture for Harish-Chandra modules. Let $\KR$ denote a
maximal compact subgroup of $\GR$ and let $K$ be the complexification of $\KR$.
Then there is an equivalence of categories:
$$
\{\text{certain $K$-equivariant sheaves on $X$}\}
@>{\ \ \sim \ \ }>>\{\text{H-C-modules}\}\,,
\tag1.2
$$
where the functor is again given by cohomology. In lectures 2 and 3 we make precise
the left hand sides of these constructions. We discuss the construction (1.1), due to
Kashiwara and Schmid \cite{KSd}, in lecture 4. The construction (1.2), due to
Beilinson-Bernstein \cite{BB1}, which is older, we only discuss briefly in lecture
5. The constructions (1.1) and (1.2) fit together as  follows:
$$
\CD
\{\text{certain $\GR$-equivariant sheaves on $X$}\} @>{\ \ \ \ \ }>>\{G_\Bbb R
\text{-representations}\}
\\
@VVV @VVV
\\
\{\text{certain $K$-equivariant sheaves on $X$}\}
@>{\ \ \sim \ \ }>>\{\text{H-C-modules}\}\,,
\endCD
\tag1.3
$$
where the second vertical arrow associates to a representation its Harish-Chandra
module. The first vertical arrow is called the Matsuki correspondence for
sheaves which we discuss in some detail in lecture 5. It was shown to be an
equivalence in \cite{MUV}.

It is clear that equivalences of categories can be used to answer
questions of categorical nature. For example, the question of how 
the standard representations decompose into irreducible representations, i.e., the
Kazhdan-Lusztig conjectures,  can be translated by equivalence (1.2) into a
question about
$K$-equivariant sheaves and was solved in this way. On the other hand, it is not
immediately clear that interesting invariants of representations can be
constructed directly from the geometric data. As an example of our techniques we
give, in lecture 7,  a geometric formula for the character of a representation.
Here it is crucial that we use the equivalence (1.1). Finally, in lecture 8, we
briefly discuss the solution of the Barbasch-Vogan conjecture, where we use the
fact that the constructions (1.1) and (1.2) fit together to form diagram (1.3).

\remark{Remark} There does not seem to be any direct way of deciding when a
representation is unitary  using the construction (1.1).  Kostant and
Kirillov have suggested that, from the point of view of unitary representations,
one should consider the dual Lie algebra $\fg_\RR^*$ instead of the flag manifold.
At this time, it is not clear what sheaves $\cf$ one should consider on
$\fg_\RR^*$. For guidance on this matter, one should consult the lectures of
Vogan. In his lectures, Vogan points out that in the case of unitary
representations and $\fg_\RR^*$ there probably  does not exist as nice a
dictionary as (1.1) and (1.2).
\endremark

To give a more detailed idea of what will be done in lectures 6-8, we first set
up some notation and then discuss the case of compact groups. 

\subheading{Notation}
We begin by introducing notation which will be used throughout this
paper. Consider the flag variety $X$ of the complexification $G$ of $\GR$,
and let us  view it as the variety of all Borel subalgebras of
$\fg=\operatorname{Lie}(G)$. The variety
$X$ carries a tautological bundle $\Cal B$ whose fiber over $x\in X$ is the Borel
subalgebra $\fb_x$ which fixes $x$. The bundle $\Cal B$ is $G$-homogenous. In
particular, it is determined by the adjoint action of the Borel group $B_x$, the
stabilizer group of $x\in X$, on the fiber $\fb_x$. From this we conclude that the
$G$-bundle $\Cal B/[\Cal B,\Cal B]$ is trivial.  Its fiber $\fh$ is called the
universal Cartan algebra; by definition it is canonically isomorphic to
$\fb_x/[\fb_x,\fb_x]$ for any $x\in X$. Any concrete Cartan
$\ft\subset \fg$ has a set of fixed points on $X$ with the same cardinality as
the Weyl group $W$. A choice of a fixed point of $x$ amounts to a choice of a
Borel $\fb_x \supset \ft$ and hence determines a canonical isomorphism $\tau_x:
\ft \to \fh$. Via these isomorphisms $\fh^*$ inherits a canonical root system
$\Phi$. Furthermore, $\Phi$ comes equipped with a canonical system of positive
roots $\Phi^+$ such that the roots of $\fg/\fb_x$ are positive.
Similarly, we can define the universal Cartan group $H\cong B_x/[B_x,B_x]$
for $G$ whose Lie algebra is $\fh$. Let $\Lambda\subset \fh^*$ denote the
$H$-integral lattice, i.e., the set of $\lambda\in\fh^*$ which lift to a
character of $H$. The elements in $\Lambda$ correspond to $G$-equivariant
holomorphic line bundles on $X$: an element $\mu\in \Lambda$ determines a
character $e^\mu$ of $H\cong B_x/[B_x,B_x]$ which lifts to $B_x$ and
hence gives rise to a $G$-equivariant line bundle $\bold L_\mu$ on $X$. We will build
the $\rho$-shift into our notation from the beginning and write $\Cal
O(\mu+\rho)$ for the sheaf of holomorphic sections of $\bold L_\mu$. Here, as
usual, $\rho\in\fh^*$ stands for half the sum of positive roots.

\subheading{The compact case}
We will now consider the case when the group $\GR$ is compact. The irreducible
representations of $\GR$ are parametrized by $\lambda\in \Lambda+\rho$ such that
$\lambda-\rho$ is dominant and the representations themselves are concretely
exhibited on the sections of the $G$-homogenous line bundle $\Cal O(\lambda)$,
i.e., on 
$$
\oh^0(X,\O(\lambda)) \ = \ \oh^0(X,\O(\bold L_{\lambda-\rho}))\ .
\tag1.4
$$
This implements constructions (1.1) and (1.2) for compact groups. 
The element $\lambda$ determines a coadjoint $\GR$-orbit 
$\Omega \subset i\fg_\RR^*$, and an isomorphism $X \cong \Omega$, as follows. Given
$x\in X$, there is a unique Cartan $\TR$ of $\GR$ which fixes $x$. As was explained
above, this gives us a map $\tau_x: \ft_\Bbb R = \operatorname{Lie}(\TR) \to
\fh$,  which, in turn, allows us to pull back $\lambda\in\fh^*$ to an 
element $\lambda_x\in i\frak t_\Bbb R^*$. The direct sum decomposition $\gr=\frak
t_\Bbb R \oplus [\frak t_\Bbb R,\gr]$ allows us to view $\lambda_x$ as an
element in $i\frak g_\Bbb R^*$. This construction provides  a $\GR$-equivariant
isomorphism between $X$ and a $\GR$-orbit $\Omega\subset i\fg_\Bbb R^*$.  As a
coadjoint orbit,  $\Omega$ has a canonical  symplectic  form $\sigma_\Omega$.
 We define the Fourier transform $\hat\phi$ of a tempered 
function $\phi$, without choosing a square root of $-1$:
$$
\hat \phi (\zeta) \ = \ \int_\gr e^{\zeta(x)}\phi(x) dx \qquad
(\zeta\in i\fg^*_\Bbb R)\,.
\tag1.5
$$
Let $\Theta$ denote the character of the representation on
$\oh^0(X,\O(\lambda))$ with highest weight $\lambda-\rho$
and let $\theta = (\det \exp_*)^{1/2}\exp^*\Theta$ denote its character on the
Lie  algebra. Then, according to Harish-Chandra,
$$
\int_{\gg_\RR}
\theta_\Omega\, \phi\, dx \ = \ \frac 1 {(2\pi i)^n n!}
\int_{\Omega}  \hat \phi \ \sigma_\Omega^n
\,.
\tag1.6
$$
In other words, 
$$
\text{Fourier transform of $\theta$}\ = \ \text{the coadjoint orbit $\Omega$
with measure $\frac{\sigma_\Omega^n}{(2\pi i)^n n!}$}\,.
\tag1.7
$$
In lecture 7 we generalize this formula for representations of an arbitrary
semisimple Lie group $\GR$. A crucial ingredient of this generalization is the
characteristic cycle construction of Kashiwara, which we discuss in lecture 6. 
In the paper \cite{SV2}, where the generalization of (1.6) is given, we also obtain
an other character formula, resembling the  Weyl's character
formula, which gives the  character of a representation via a Lefschetz type fixed
point formula.

Our first goal, in lectures 2-4, is to define the functors
$$
\{\text{$\GR$-equivariant sheaves on $X$}\} \longrightarrow\{G_\Bbb R
\text{-representations}\}
\tag1.8
$$
To do so, we first must give a precise meaning to the left hand side.
It is the ``twisted'' $\GR$-equivariant derived category 
$D_\GR(X)_\lambda$ of constructible sheaves on the flag manifold $X$. Here the
twisting parameter
$\lambda
\in
\frak
\fh^* =$ the dual space of the universal Cartan. We will explain the construction of
$D_\GR(X)_\lambda$ in three stages. 
First we introduce the notion of the derived category of constructible sheaves,
then we describe how to make this notion $\GR$-equivariant, and finally, we 
explain the  twisted version. Because the $\GR$-orbits on the flag manifold are
semi-algebraic sets, we will develop the general theory in this context.


\Lecture2{ Derived categories of constructible sheaves}			%


In this lecture we give a brief treatment of constructible sheaves and derived
categories. At the end of this lecture series there is an appendix by Markus Hunziker
which one may consult for further basic information about derived categories and sheaf
cohomology.  For a more detailed discussion see, for example,
\cite{KSa}.

\subhead Semi-algebraic sets \endsubhead
Recall that a subset of $\RR^n$ is called semi-algebraic if 
it is the union of finitely many sets of the form 
$$
  S\ = \ \{x \in \RR^n \mid f_1(x)=\cdots = f_r(x)=0\ ,\ g_1(x) > 0 \ , \ldots ,\
g_s(x)>0\, \},
$$
where the $f_i$ and $g_j$ are polynomials in $\RR[x_1,\ldots,x_n]$.
It follows directly from the definition that the class of semi-algebraic sets
is stable under finite intersections, finite unions, and taking the complement.
If $S \subset \RR^n$ and $S' \subset \RR^m$ are semi-algebraic sets then
a map $f: S \rightarrow S'$ is called semi-algebraic if it is continuous
and its graph is a semi-algebraic set in $\RR^{n+m}$. The composition of two
semi-algebraic maps  is semialgebraic, of course. Furthermore: 
$$
\aligned\text{the image of a semi-algebraic
set under}
\\
\text{a semi-algebraic map is  semi-algebraic.}
\endaligned
\tag2.1
$$ 
For this fact see, for example, \cite{BM}.
The definition of a semi-algebraic set generalizes readily to arbitrary algebraic 
manifolds as follows.
If $X$ is a real algebraic manifold then a subset $S$ of $X$ is called semi-algebraic if
$S\cap U$ is semi-algebraic for every Zariski open affine subset $U$ of $X$. 
The following crucial, non-trivial property of semi-algebraic sets will be important to 
us:
$$
 \text{every semi-algebraic set can be triangulated by semi-algebraic simplices.}
\tag2.2
$$ 
For a proof of this fact in a far more general context, see \cite{DM, 4.10}.
\remark{2.3 Remark} In Lecture 8 we have to work in a more general setting of geometric
categories arising from the the o-minimal structure $\Bbb R_{\text{an,exp}}$. The
article \cite{DM} provides an excellent exposition of this theory. 
\endremark

\subhead Constructible sheaves \endsubhead

Let $X$ be a semi-algebraic set and let us fix a semi-algebraic triangulation $\T$ of
$X$. A sheaf $\cf$ on
$X$ is called $\T$-constructible if its restriction to  any (open) simplex $\sigma$
of $\T$  is a constant sheaf of complex vector spaces of finite rank. We call a sheaf
constructible if it is $\T$-constructible for some $\T$. 
The notion of a  $\T$-constructible sheaf is the natural notion of a
coefficient system on a simplicial complex: it associates a vector space to each 
(open) simplex and a linear map from  the vector space of one simplex to that of 
another if the first simplex lies on the boundary of the second.

Another way of thinking about constructible sheaves is based on the notion of
a local system. Recall that   local systems  are locally constant sheaves of
finite rank. They constitute a special case of constructible sheaves. Given a
locally  constant sheaf $\cf$ on $X$ and a point $x\in X$ we obtain a representation
of
$\pi_1(X,x)$ on the stalk $\cf_x$ by continuing the sections of $\F_x$ along the loops
based at $x$. This gives (for $X$ connected!) an equivalence of categories:
$$
  \{\text{local systems on $X$}\} \leftrightarrow
  \{\text{finite dimensional representations of $\pi_1(X,x)$}\}
\tag2.4
$$
The category of constructible sheaves on $X$ is the smallest abelian
category containing the local systems on all semi-algebraic subsets of $X$.

\subhead Derived categories \endsubhead

Derived categories provide a convenient tool which helps to organize arguments in
homological algebra. They provide the appropriate framework for resolutions and for
derived functors. The fundamental idea is that one works systematically, from the
outset, in the category of complexes. Also, there are functors that exist only
in the context of derived categories. An example is provided by the functor $f^!$
introduced later in this section. 

Let  $X$ be an arbitrary semi-algebraic
set. The bounded derived category  $\od(X)$ of constructible sheaves on $X$ has as
its objects bounded complexes of constructible sheaves. Its morphisms are given by
chain homotopy classes of maps of chain complexes and, in addition, we formally invert
the maps
$\phi:\Cal A^\bullet \to \Cal B^\bullet$ (the quasi-isomorphisms) which induce 
isomorphisms $\oh^\bullet(\phi):\oh^\bullet(\Cal B^\bullet) @>{\ \sim\ }>>
\oh^\bullet(\Cal A^\bullet)$  on the cohomology sheaves. The notion of exactness
looses its meaning in $\od(X)$. Exact sequences of chain complexes are called
distinguished triangles when viewed in $\od(X)$. Objects in $\od(X)$ can be shifted: if
$\Cal A^\bullet\in\od(X)$ then $\Cal A^\bullet[n]$ denotes the complex such
that\footnote{One should also multiply the differential of the complex $\Cal A^\bullet$
by $(-1)^{n}$}
$(\Cal A^\bullet[n])^k=\Cal A^{n+k}$.

Let us consider complexes of $\T$-constructible sheaves and let us denote the resulting
derived category by $\od_\T(X)$. The injective $\T$-constructible sheaves   are easy to
describe. Any injective $\T$-constructible sheaf is a direct sums of basic injective
$\Cal{T}$-constructible sheaves. A basic
injective
$\Cal{T}$-constructible  sheaf is a constant sheaf on a closure of a  simplex in
$\Cal{T}$. To see that a constant sheaf on closure of simplex $\sigma$ is injective,
it suffices to note that
$$
\Hom(\bc_{\bar\sigma}, \cf) \ = \ (\cf_\sigma)^*\,, \qquad\text{for any $\Cal
T$-constructible sheaf $\cf$}\,.
$$
From the discussion above, it follows easily that the category of
$\T$-constructible sheaves has enough injectives. In particular, every $\F \in
\od_\T(X)$ is isomorphic (in
$\od_\T(X)$) to a complex of injectives, its injective resolution. Furthermore, by
a standard argument (see, for example, \cite{KSa, Prop. 1.8.7}), 
$$
  \od_\T(X) \ \cong \ \text{homotopy category of injective $\T$-complexes}\,.
\tag2.5
$$
Except for trivial cases, the category of constructible sheaves on $X$ does not have
enough injectives (exercise). However, any $\cf\in\od(X)$ lies in some $\od_\T(X)$  and
hence has an injective representative in $\od_\T(X)$. One can develop the theory for
various operations on sheaves utilizing  this principle. However, it is technically 
simpler, and perhaps more elegant, to view $\od(X)$ as a subcategory of the
derived category of all sheaves of $\Bbb C$-vector spaces on $X$ and take the injective
representative of $\cf\in\od(X)$ inside this bigger derived category. In what follows,
we will take this point of view. In particular, we view $\od(X)$ as a subcategory of
the derived category of all sheaves of $\Bbb C$-vector spaces on $X$ consisting of 
complexes $\cf$ such that  the cohomology sheaves $\oh^k(\cf)$ are
constructible and are non-zero for finitely many values of $k$ only.

\subhead Operations on sheaves \endsubhead

From now on  all of the semi-algebraic sets are assumed to be locally compact. Let
$f: X
\rightarrow Y$ be a map of semi-algebraic sets.  We shall define functors $Rf_*,Rf_! :
\od(X)
\rightarrow
\od(Y)$ and 
$f^*, f^! : \od(Y) \rightarrow \od(X)$.

\definition{Direct image} 
If $\F$ is any sheaf on $X$ then the direct image of $\F$ by $f$, denoted
by $f_*\F$ is the sheaf on $Y$ defined by:
$$
  V \longmapsto f_*\F(V) := \F(f^{-1}(V))
$$
It is not immediately clear that the constructibility of $\F$ implies the
constructibility of 
$f_*(\F)$. It follows from that fact that any semi-algebraic map can
be Whitney stratified\footnote{Loosely speaking this means that the space Y can be
decomposed into strata  so that the map $f$ restricted to
$f^{-1}(S)$ is locally trivial for any  stratum $S$ of $Y$.}. For a very general
discussion of such matters, see
\cite{DM}. The functor
$f_*$ lifts to a functor
$Rf_*:
\od(X)
\rightarrow
\od(Y)$ in the usual way: if $\Cal J$ is an injective
resolution of an object $\F\in \od(X)$ then 
$Rf_*(\F) = f_* \Cal J$. Here, again, the fact that $Rf_*(\F)$ is constructible
is not entirely obvious; recall that the injective resolution $\Cal J$ is not
a complex of constructible sheaves. To see that $Rf_*(\F)\in \od(X)$, one can argue in
the same way as proving the constructibility of $f_*\cf$.        If
$f: X
\rightarrow
\{\pt\}$ then
$f_* =
\Gamma(X, \ \,)$ is the global section functor on sheaves and hence the complex 
$Rf_*(\CC_X) = R\Gamma(X,\CC_X)$ computes the cohomology of $X$
(with coefficients in $\CC$):
$$
  \oh^k(X,\CC_X)\  =\  R^k\Gamma(X,\CC_X) \ =\  R^kf_*(\CC_X)\,.
\tag2.7
$$
The pushforward construction is functorial in the sense that
$$
R(f\circ g)_* \ = \ Rf_*\circ Rg_* \qquad \text{when} \qquad X@>f>> Y @>g>> Z\,.
\tag2.8
$$
\xca{Exercise}
Assume that $X$ is compact manifold and $f: X \rightarrow \RR$\ \, is a Morse function. 
Describe the complex $Rf_*\CC_X$. 
\endxca

\enddefinition

\definition{Inverse image}
The direct image functor $f_*$ has a left adjoint functor $f^*$ in the 
category of (constructible) sheaves, i.e., 
for any sheaf $\F$ on $X$ and any sheaf $\G$ on $Y$ one has
$$
\Hom(f^*\G,\F)\ = \ \Hom(\G,f_*\F)
\tag2.9
$$
 The sheaf $f^*\G$ is called the inverse 
image of $\G$ by $f$. An explicit construction of $f^*\G$ 
is as follows: $f^*\G$ is the sheaf  associated to the presheaf 
$$
U \longmapsto \ \lim\Sb   @>>> \\ V\supset f(U) \endSb \G(V) \ .
$$ 
It is clear from this description of $f^*\G$ that if $x\in X$ then  
$$
  (f^{*}\G)_x\  = \ \G_{f(x)}\ .
$$ 
This implies that $f^*$ is exact and hence $Rf^*\cf=f^*\cf$, for any $\cf\in\od(X)$.
\enddefinition

\definition{Direct image with proper support}
For a (constructible)
sheaf  $\F$  on $X$  we define its direct image with proper support $f_!\F$ as the
following subsheaf of $f_*\F$:
$$
  V \longmapsto\  f_!\F(V)\  =\  \{s \in \F(f^{-1}(V)) \mid f:\supp(s) \rightarrow U\  
   		                        \text{is proper} \}
$$
Note that if the map $f$ is proper then $f_*\F = f_!\F$. In the other extreme, if $j:U
\to X$ an  embedding then 
$$
j_!\cf \ = \ \text{extension of $\cf$ by zero}\,.
\tag2.10
$$
The functor $f_!$ is closely related to cohomology with compact support.
Let 
$$
  \Gamma_c(X,\F)\  = \ \text{global sections of $\F$ with compact support}
$$
If $f:X \rightarrow \{\pt\}$ then $f_!\cf =\Gamma_c(X,\F)$ and hence
$$
\oh^k_c(X,\CC_X)\  =\  R^k\Gamma_c(X,\CC_X) \ =\  R^k\!f_!(\CC_X)\,.
\tag2.11
$$
\enddefinition

\definition{Base change}
Consider a Cartesian square of semi-algebraic sets:
$$
  \CD
  X' @>u>> X \\
  @V{v}V{\qquad \qed}V  @VV{f}V \\
  Y'  @>g>> Y \,.
 \endCD
$$
Recall that a square is Cartesian if it commutes and  $X' \simeq
X\times_Y Y'$. Then  we have a natural isomorphism of functors:
$$
  g^*\circ Rf_! \ \simeq\  Rv_! \circ u^*
\tag2.12
$$
In particular, if $Y'=\{y\}$, then (2.12) implies that
$$
  (R^kf_!\F)_y \ \simeq \ R^k\Gamma_c(f^{-1}(y), \F|_{f^{-1}(y)})\ \simeq \
\oh^k_c(f^{-1}(y),\cf) .
$$
\enddefinition

\head Verdier duality \endhead

\noindent
Unlike $f_*$, the functor $f_!$ does not have a right adjoint within the
category of sheaves. However, the functor $Rf_!: \od(X) \rightarrow
\od(Y)$ does have a right adjoint functor $f^!: \od(Y) \rightarrow \od(X)$:
$$
  \RHom(Rf_!\F,\G) \ = \ \RHom(\F,f^!\G)\, .
\tag2.13
$$
This statement is usually referred to as Verdier duality.
The functor $f^!$ can not be obtained by taking the derived functor of a functor on
sheaves. If
$i:Y
\hookrightarrow X$ is locally closed embedding then
$$
i^!\cf \ = \ R\Gamma_Y(\cf)\,,
$$
where 
$$
\Gamma_Y(\cf) \ = \ \text{sections of $\cf$ supported on $Y$}\,.
$$
In particular,
$$
\oh^k(X,i^!\cf)\ = \ \oh^k(X,R\Gamma_Y(\cf))\ = \ \oh^k_Y(X,\cf)\ = \ \oh^k(X,X-Y;\cf)
$$
the local cohomology of $\cf$ along $Y$.

\definition{Dualizing complex and duality functor} 
Let $X$ be a semi-algebraic set, and $f: X \rightarrow \{\pt\}$.
We define the dualizing complex as
$$
{\Bbb D}_X\  =_{\text{def}} \ f^! \CC_{\{\pt\}}\,.
\tag2.14
$$
Applying Verdier duality to the map $f: X \rightarrow \{\pt\}$, we get
$$
\RHom(Rf_!\Bbb C_X, \CC_{\{\pt\}}) \ \cong \ \RHom(\Bbb C_X,{\Bbb D}_X) \ \cong \
R\Gamma(X,{\Bbb D}_X)\,.
$$
In particular, we see that 
$$
\oh^k(X,{\Bbb D}_X)\ \cong \  (\oh^{-k}_c(X,\CC_X))^* \ \cong\ \oh_{-k}(X,\Bbb C_X)\,.
$$
By a slightly more refined computation we see that the dualizing complex can be
interpreted as the complex of chains -- with closed, not necessarily compact, supports
-- on
$X$, placed in negative degrees. In particular, there is a canonical isomorphism
$$
\oh^k(\Bbb D_X)_x \ \cong \ \oh_{-k}(X,X-\{x\};\Bbb C)\,, \qquad \text{for every
$x\in X$}
\tag2.15
$$
For $\F \in \od(X)$ we define  $\Dual_X \F$, the Verdier dual of $\F$, by the formula
$$
  \Dual_X \F\  =_{\text{def}} \ \RShHom(\F,{\Bbb D}_X)\,.
\tag2.16
$$
Here the functor $\Cal Hom$ is
the sheafification of the functor $U\mapsto \Hom(\cf|U,\cg|U)$. If $\F \in \od(X)$
then $\Dual_X \F \in \od(X)$, and 
$$ 
\F @>{\ \sim \ }>> \Dual_X\Dual_X \F
\tag2.17
$$
is an isomorphism. This latter statement is called biduality. Furthermore, if
$f:X\rightarrow Y$ is a morphism of
semi-algebraic sets, then:
$$
  Rf_! = \Dual_X\circ Rf_* \circ \Dual_Y  \qquad \text{and}
  \qquad 
  f^!  = \Dual_Y\circ f^*  \circ \Dual_X \ .
\tag2.18
$$
\enddefinition

\definition{Poincar\'e duality}
Let $X$ be an oriented (semi-algebraic) manifold of dimension $n$, 
and let $f: X \rightarrow \{\pt\}$.
Then, by  formula (2.15), $\oh^k(f^!\CC_{\{\pt\}})\cong \oh_{-k}(X,X-\{x\};\Bbb C)$.
Hence, the specific orientation of $X$ provides a distinguished isomorphism
$$
f^!\CC_{\{\pt\}} = \Bbb D_X \cong\CC_X[n]\,.
\tag2.19
$$
By Verdier duality
$$
  \RHom(\RGamma_c(X,\CC_X)[n],\CC) \simeq \RGamma(X,\CC_X)\,,
\tag2.20
$$
and taking the $p$-th cohomology group gives the isomorphism:
$$
  (\oh^{n-p}_c(X,\CC_X))^*\simeq \oh^p(X,\CC_X)\,.
\tag2.21
$$
\enddefinition

\Lecture3{Equivariant derived categories}			

Let $X$ be a semi-algebraic set with an algebraic $\GR$-action.
Bernstein and Lunts \cite{BL} defined a category $\od_\GR(X)$, 
the $\GR$-equivariant derived category of  constructible sheaves on $X$, 
together with a forgetful functor $\od_\GR(X) \rightarrow \od(X)$, satisfying
the following condition:
if $f: X \rightarrow Y$ is a $\GR$-equivariant map between 
semi-algebraic $\GR$-spaces, then the functors 
$Rf_*,Rf_!,f^*,f^!$ lift canonically to functors between the $\GR$-equivariant
derived categories. Here we give a brief account of this theory. For a short summary,
see also \cite{MV}.

\head $\GR$-equivariant sheaves \endhead

As before, we  assume that $\GR$ is a connected semisimple Lie group and let us
consider a semi-algebraic $\GR$-space $X$. Naively, a
$\GR$-equivariant (constructible) sheaf on
$X$ is a sheaf
$\F$ on
$X$ together with isomorphisms of the stalks
$$
  \phi_{(g,x)}\ : \ \cf_{gx} @>{\ \sim\ }>> \cf_x \qquad \text{for all} \
g\in\GR,\, x\in X\,.
\tag3.1
$$ 
Of course, one wants the isomorphisms $\phi_{(g,x)}$ to depend continuously on $g$ and
$x$. To make the notion of a $\GR$-equivariant map precise consider the maps 
$a,p: \GR\times X \rightarrow X$, $a(g,x)=gx$, $p(g,x)=x$.
A $\GR$-equivariant sheaf on $X$ is sheaf $\F$ on $X$ together with an isomorphism
of sheaves on $\GR\times X$\,:
$$
\phi\ :\  a^*\F \ @>{\ \sim  \ }>> \ p^*\F
\tag3.2
$$
such that $\phi|_{\{e\}\times X} = \text{id}$. (If $\GR$ is not assumed to be connected
we  need to add a cocycle condition.) Clearly, by restricting to the stalk at $(g,x)$,
(3.2) gives (3.1).

\xca{Exercise}
Let $\F$ be a $\GR$-equivariant sheaf on $X$. Construct a canonical linear
$\GR$-action  on the cohomology spaces $\oh^k(X,\F)$.
\endxca

\subheading{$\GR$-equivariant local systems}
If $\GR$ acts transitively on $X$ then, as is not difficult to see,  we have an
equivalence of categories
$$
\left\{\matrix\text{$\GR$-equivariant}\\
 \text{sheaves on $X$}\endmatrix\right\}\ \longleftrightarrow\ 
 \left\{\matrix\text{finite dimensional representations of}\\
        \text{the component group $(\GR)_x/(\GR)_x^0$}\endmatrix\right\}\,.
\tag3.3
$$
Fix $x\in X$. Then we have a fibration $\GR \to X$, $g\mapsto gx$, with fiber the
stabilizer $(\GR)_x$. Because $\GR$ is connected, the long exact sequence of homotopy
groups for a fibration yields a surjection  $\pi_1(X,x)\to\pi_0((\GR)_x) =
(\GR)_x/(\GR)_x^0$. Thus, a $\GR$-equivariant sheaf gives rise to a local system
on $X$.

\head Equivariant derived categories and functors  \endhead

Let $X$ be a semi-algebraic  $\GR$-space. The equivariant derived category can be
characterized by the following properties:

\roster
\item "(1)"There exists a forgetful functor $\od_\GR(X) \rightarrow \od(X)$
which maps $\GR$-equivariant sheaves to constructible sheaves.

\item "(2)"For a $\GR$-equivariant map $f: X \rightarrow Y$ the functors
$Rf_*,Rf_!. f^*,f^!$ lift to functors between $\od_\GR(X)$ and $\od_\GR(Y)$.
The same is true for the duality functor $\Dual$ and the standard
properties  hold for these lifted functors.

\item "(3)"
If $G'_\RR\subset \GR$ is a normal subgroup and $G'_\RR$ acts freely on $X$ then 
$$\od_\GR(X)  \ \simeq \ \od_{\GR/G'_\RR}(G'_\RR\backslash X)\ .$$
\endroster

In the above, by a ``standard property" we mean properties like (2.9), (2.12), (2.13),
(2.16), and (2.17).  The idea of the construction of
$\od_\GR(X)$ is the same as that of equivariant cohomology:
if the action of $\GR$ on $X$ is free then, by property (3),  we simply set
$\od_\GR(X) = \od(\GR\backslash X)$. If the action of $\GR$ on $X$ is not free, we
replace $X$ by $X \times EG_\RR$  where $EG_\RR$ is a contractible (infinite dimensional) space on which 
$\GR$ acts freely and consider the diagonal action of $\GR$ on $X \times EG_\RR$. Let us
organize the various spaces in the following diagram: 
$$
\CD
X @<p<< X \times E\GR @>q>> \GR\backslash(X\times E\GR)\,,
\endCD
\tag3.4
$$
where $p$ is the projection to the first factor and $q$ is the map to the quotient. We
can now make the following formal definition.  An object in $\od_\GR(X)$ is   a
triple
$(\F,\G,\psi)$, where $\F\in \od(X)$, $\Cal{G}\in\od(\GR\backslash(X\times E\GR))$, and
$\psi$ is an 
isomorphism, 
$$ 
\psi\,:\,p^* \F @>{\ \sim\ }>> q^* \Cal{G}\ .
\tag3.5
$$
Here $p: X \times E\GR \rightarrow X$ is the projection on the first 
factor and 
$q: X \times E\GR \rightarrow \GR\backslash(X\times E\GR)$ is the quotient 
map. To make this work in our semi-algebraic context, one approximates $E\GR$ by finite
dimensional spaces. For details, see \cite{BL} and also \cite{MV}.

\xca{Exercise}
Let $\F\in\od_\GR(X)$ and $a,p:\GR\times X \rightarrow X$ the action map and the
projection. Show that the above definition of $\od_\GR(X)$ yields a morphism $\phi:
a^*\F \rightarrow p^*\F$  in $\od(\GR\times X)$. Show that this gives  a linear
$\GR$-action on $\oh^*(X,\F)$; here $\cf$ is viewed non-equivariantly, i.e., as an
element of $\od(X)$. 
\endxca

\head Twisting \endhead

Let us return to the situation of lecture 1 and the notation used there. Now $X$ will
stand for the flag manifold on the complexification $G$ of $\GR$. The enhanced flag
variety $\hat X$ is defined as
$$
\h X \ = \ G/N \,, \qquad \text{where}\ \ N=\text{unipotent radical of a Borel $B$}
\tag3.6
$$
The group $G\times H$, where $H\cong B/N$ denotes the universal Cartan group, acts
transitively on
$\h X$ by the formula $(g,h)\cdot g'N = gg'h^{-1}N$.  A ``sheaf
with twist $\lambda \in \fh^*$ on $X$'' is a sheaf $\F$ on  $\hat{X}$ such that for any
$\hat{x}\in\hat{X}$ the pullback of $\F$ to $H$ under $h\mapsto h\cdot \hat{x}$ is
locally constant and has the same  monodromy as the function $e^{\lambda-\rho}$. We will
think of twisted sheaves as objects on $X$. Note that if
$\lambda=\rho$ or, more generally, if $\lambda$ is an
$H$-integral translate of $\rho$, then the $\lambda$-twisted sheaves are just
ordinary sheaves on $X$. The notions of lecture 2, as well as the notions of the
equivariant derived category, extend readily to the twisted case. In particular, we
have the notion of the
$\lambda$-twisted,
$\GR$-equivariant derived category $\od_\GR(X)_\lambda$. If $\mu \in \Lambda$, i.e., if
$\mu$ is $H$-integral, then $\od_\GR(X)_\lambda = \od_\GR(X)_{\lambda+\mu}$. The derived
category  $\od_\GR(X)_{\lambda}$ is generated\footnote{In the sense of a triangulated
category, i.e., by shifting and forming distinguished triangles.} by standard
sheaves. For technical reasons, which will become apparent in the next lecture, we give
a classification of standard sheaves in $\od_\GR(X)_{-\lambda}$. By definition,
standard sheaves are associated to pairs $(X,\Cal L)$, where $S$ is a $\GR$-orbit on
$X$ and $\Cal L$ is an irreducible, $(-\lambda)$-twisted $\GR$-equivariant local system
on $S$. Given such a pair $(S,\Cal L)$ we can attach to it two types of standard
sheaves:
$$
Rj_*\Cal L \ \ \ \text{and}\ \ \ \ j_!\Cal L, \qquad \ \text{where $j:S\hookrightarrow
X$ denotes the inclusion}\,.
\tag3.7
$$
We can use either type to generate $\od_\GR(X)_{-\lambda}$.
Let us fix $x\in S$ and a Cartan $\TR\subset \GR$ which fixes $x$. As was shown in
lecture 1, this data gives  an identification $\tau_x : \ft
@>{\sim}>> \fh$, where $\ft$ denotes the complexification of the Lie algebra of $\TR$.
Then
$$
\left\{\matrix\text{Irreducible, $\GR$-equivariant}\\
 \text{$(-\lambda)$-twisted local systems on $S$}\endmatrix\right\}\
\longleftrightarrow\ 
 \left\{\matrix\text{characters $\chi: \TR \to \Bbb C^*$}\\
       \text{with }\  d\chi = \tau_x^*(\lambda-\rho)\endmatrix\right\}\,.
\tag3.8
$$
To verify this statement, we  apply (3.3) for the action of $\GR\times\fh$ on $\hat S$.
Here $\hat S$ stands for the inverse image of $S$ in $\h X$ and $\fh$ acts on $\hat X$
via the exponential map $\exp:\fh\to H$. Statement (3.3) can be applied because 
$(-\lambda)$-twisted, $\GR$-equivariant local systems on $S$ are
$\GR\times\fh$-equivariant local systems on $\hat S$.

\Lecture4{Functors to representations}			

In this lecture we make precise the ideas of lecture 1 and define the functors from
our geometric parameter space $\od_\GR(X)_{-\lambda}$ to representations. Recall that
by a representation of a semisimple Lie group $\GR$ we mean a representation on a
complete, locally convex, Hausdorff topological vector space of finite length. We
denote the infinitesimal character corresponding to
$\lambda\in\fh^*$ by $\chi_\lambda$ and use Harish-Chandra's normalization so that
the value $\lambda=\rho$ corresponds to the trivial infinitesimal character. 

Let $\O(\lambda)$ denote the sheaf of $\lambda$-twisted holomorphic
functions on $X$. More precisely,  $\O(\lambda)$ is the subsheaf of $\O_{\hat{X}}$
which consists of functions whose restriction to any fiber of the map $\h X \to X$ is a
constant multiple of the function $e^{\lambda-\rho}$ when we identify the fibers of 
$\h X \to X$  with $H$ via the map $h\mapsto h\cdot gN = gh^{-1}N$, as usual. We note
that if 
$\lambda \in \Lambda +\rho$, the $\rho$-translate of the $H$-integral lattice, then
$\O(\lambda)$ can be viewed as an ordinary sheaf on $X$ and it coincides with the
sheaf of sections of the line bundle $\bold L_{\lambda-\rho}$.

In \cite{KSd}, Kashiwara and Schmid define  two functors $M$ and $m$
$$
\od_\GR(X)_{-\lambda} \ \longrightarrow\ \left\{\matrix\text{virtual admissible
$\GR$-representations of}\\ \text{finite length with infinitesimal character
$\chi_\lambda$}\endmatrix \right\}
\tag4.1a
$$
given by the formulas
$$
\aligned
 &M\ :\  \F \longmapsto \sum(-1)^k \Ext^k(\Dual\cf,\O(\lambda)) 
\\
&m\ :\  \F \longmapsto \sum(-1)^k \oh^k(X, \cf\otimes\O(\lambda)) \,.
\endaligned
\tag4.1b
$$

A few comments are in order. In \cite{KSd} the target of the functors (4.1) is a
derived category of representations. Here we have introduced the simplified
version only where the functors land in virtual representations, as we do not need the
more refined version in these lectures. To explain the formulas (4.1b), note that the
duality functor  $\Dual$ takes $\od_\GR(X)_{-\lambda}$ to $\od_\GR(X)_{\lambda}$.
The  groups $\Ext^k(\Dual\cf,\O(\lambda))$ are to be interpreted as $\Ext$-groups in the
category of
$\lambda$-twisted sheaves. As to the groups $\oh^k(X, \cf\otimes\O(\lambda))$, note
that the sheaf $\cf\otimes\O(\lambda)$ has trivial monodromy along the fibers of
$\hat X\to X$. Therefore  $\cf\otimes\O(\lambda)$ descends to $X$ and $\oh^k(X,
\cf\otimes\O(\lambda))$ stands simply for cohomology on $X$. The topology on the
spaces  $\Ext^k(\Dual\cf,\O(\lambda))$ and $\oh^k(X, \cf\otimes\O(\lambda))$ is induced
by the usual topology on the sheaf of holomorphic functions $\O(\lambda)$. 
Although in general the formulas (4.1) define 
functors into virtual representations only, there is a subcategory of
$\od_\GR(X)_{-\lambda}$ such that restricted to this subcategory the groups
$\Ext^k(\Dual\cf,\O(\lambda))$ and $\oh^k(X,
\cf\otimes\O(\lambda))$ are non-zero only in degree zero. Restricted to this
subcategory, the functors $M$ and $m$ land in representations, and every
representation, up to infinitesimal equivalence, arises in this fashion. We will
discuss this subcategory in more detail in lecture 5. Moreover, given $\cf$, the
representations $M(\cf)$ and
$m(\cf)$ are the maximal and minimal globalizations of Schmid \cite{S1} of the same
Harish-Chandra module. Strictly speaking, this statement is not made in \cite{KSd}.
However, it follows from the results in \cite{KSd} and the statement (4.2) below. 
As was noted earlier, the categories $\od_\GR(X)_{-\lambda}$
and  $\od_\GR(X)_{-\lambda -\mu}$ are canonically equivalent if $\mu$ lies in
the $H$-integral lattice $\Lambda$. Passing from the parameter value $\lambda$ to
$\lambda+\mu$ in the construction (4.1) amounts to coherent continuation on the
representation theoretic side. 

The functors $M$ and $m$ turn the duality operation $\Bbb D$ into duality of
representations:
$$
\gathered
\text{The virtual representation $\sum
\,(-1)^p\,\operatorname{Ext}^p(\Bbb D\Cal F, \Cal O_X(\lambda))$ is, up
to }
\\
\text{infinitesimal equivalence, the dual of}\  \sum
\,(-1)^{n+p}\,\operatorname{Ext}^{p}(\Cal F, \Cal O_X(-\lambda))\,.
\endgathered
\tag4.2
$$
In the above statement, $n$ stands for the complex dimension of $X$. This statement
does not appear in \cite{KSd} but a short argument for it can be found in \cite{SV2}.

The simplest example of the functors (4.1) is the case when $\cf$ is the
constant sheaf $\Bbb C_X$. For the constant sheaf to be $(-\lambda)$-twisted, the
parameter $\lambda$ has to lie in $\Lambda+\rho$. In this case, the functors (4.1)
 coincide and amount to the Borel-Weil-Bott realization of a finite
dimensional representation of $\GR$.

Next we apply the functors $M$ and $m$ to standard sheaves. To that end, let $\Cal L$
be an irreducible $\GR$-equivariant local system with twist $-\lambda-\rho$ on a
$\GR$-orbit $S$, let $j: S \hookrightarrow X$ denote the inclusion, and set $\F=
Rj_*\Cal L$.  Let us first deal with the case of discrete series and hence suppose
that $\GR$ has a compact Cartan. Furthermore, we assume that the orbit $S$ is open
in $X$. By (3.8) there is a $\GR$-equivariant $(-\lambda-\rho)$-twisted irreducible
local system on $S$ precisely when $\lambda\in\Lambda+\rho$, and such a local system is
, by necessity, the  trivial one. Hence, in our special case, 
 $\cf=Rj_*\Bbb C_S$. The representation associated to $\cf$ by (4.1) is a
discrete series representation if $\lambda$ is regular antidominant, as the following
calculation will show. This coincides with the usual  parametrization of the discrete
series representations as explained in the lectures of Zierau. To do the calculation, we
first note that, as
$\GR$-equivariant sheaves, 
$$
  \Dual \F \cong \Dual Rj_*\Bbb C_S \cong j_!\Dual_S \cong  j_!\Bbb C_S[2n]\,,
\tag4.3
$$ 
where $n=\dim_\Bbb C X=\dim_\Bbb C S$. From (4.3) and using (2.13) we conclude:
$$
\aligned 
\Ext^k(\Dual\F,\O(\lambda)) &= \Ext^k(j_!\Bbb C_S[2n],
\O(\lambda)) \\ &= \Ext^{k-2n}(j_!\Bbb C_S,\O(\lambda)) \\
&= \Ext^{k-2n}(\CC_S, 
\O(\lambda)|S)
\\
&= \oh^{k-2n}(S,\O(\lambda))
\,;
\endaligned
\tag4.4
$$
The cohomology groups $\oh^{k-2n}(S,\O(\lambda))$ are non-zero precisely when
$s=k-2n=\frac 1 2 \dim\KR/\TR$, and give a discrete series representation of $\GR$. For
this fact, see the lectures of Zierau. We thus get
$$
M(\cf) \ = \ (-1)^s\{\text{discrete series representation attached to
$(S,\lambda)$}\}\,,
\tag4.5
$$
as virtual representations. Let us now return to the case of a general $\GR$-orbit
$S$ with a  $(-\lambda-\rho)$-twisted, $\GR$-equivariant
local system $\Cal L$  on $S$, and we set  $\cf= Rj_*\Cal L$. Then,
$$
  \Dual \F = \Dual Rj_*\Cal L = j_!\Dual \Cal L = j_!\Cal L^*\otimes\text{or}_S[\dim
S]\,,
\tag4.6
$$
where $\Cal L^* =
\Hom(\Cal L, \Bbb C_S)$ denotes the dual of $\Cal L$ as a local system and
$\text{or}_S$ denotes the orientation sheaf of $S$. From (4.6) and using (2.13) we
conclude:
$$
\aligned 
\Ext^k(\Dual\F,\O(\lambda)) &= \Ext^k(j_!(\Cal L^*\otimes\text{or}_S[\dim S]),
\O(\lambda)) \\ &= \Ext^{k-\dim S}(j_! (\Cal L^*\otimes\text{or}_S),\O(\lambda)) \\
&= \Ext^{k-\dim S}(\Cal L^*\otimes\text{or}_S,j^!\O(\lambda))  \\
&= \Ext^{k-\dim S}(\CC_S,\Cal L\otimes\text{or}_S\otimes j^! \O(\lambda))  \\
&= \Ext^{k-\dim S}(\CC_S,j^!(\tilde\Cal L\otimes\widetilde{\text{or}_S}\otimes 
\O(\lambda)))
\\
&= \oh^{k-\dim S}_S(X,\tilde\Cal L\otimes\widetilde{\text{or}_S}\otimes\O(\lambda))
\,.
\endaligned
\tag4.7
$$
Here $\tilde\Cal L$ and $\widetilde{\text{or}_S}$ denote extensions of the sheaves
$\Cal L$ and ${\text{or}_S}$ to a small neighborhood of $S$. For $\lambda$
antidominant, these groups are non-zero in one degree only \cite{SW}.

Let us turn to the functor $m$. Attempting to apply $m$ to the standard sheaf
$Rj_*\cf$ leads to a seemingly very difficult calculation. However it is easy to
apply it to the standard sheaf $\cf=j_!\Cal L$. This gives
$$
\aligned 
\oh^k(X,\F\otimes \O(\lambda)) &= \oh^k(X,j_!\Cal L\otimes \O(\lambda)) \\ &=
\oh^k(X,j_!(\Cal L\otimes \O(\lambda)|S)) \\ &= \oh^k_c(X,j_!(\Cal L\otimes
\O(\lambda)|S))  \\ &= \oh^k_c(S,\Cal L\otimes
\O(\lambda)|S) \,.
\endaligned
\tag4.8
$$
If we take, in (4.8), $\cf= \Bbb DRj_*\Cal L = j_!(\Cal
L^*\otimes{\text{or}_S})[\dim S]\in \od_\GR(X)_\lambda$, then we get
$$
\oh^k(X,\F\otimes \O(-\lambda)) \ = \ \oh^{k-\dim S}_c(S,\Cal
L^*\otimes{\text{or}_S}\otimes
\O(-\lambda)|S) \,,
\tag4.9
$$
which is, by (4.2), dual to the representation in (4.7). Note that
there is a pairing
$$
(\Cal
L\otimes{\text{or}_S}\otimes\O(\lambda))\otimes(\Cal L^*\otimes{\text{or}_S}\otimes
\O(-\lambda)|S)) \longrightarrow \bold L_{-2\rho}=\Omega_X\,.
\tag4.10
$$
Hence, the duality between (4.7) and (4.9) is an extension of Serre duality. 

In (4.1) the categories $\od_\GR(X)_{-\lambda}$ and $\od_\GR(X)_{-\mu}$ map to
representations of the same infinitesimal character if $\mu$ lies in the $W$-orbit
of $\lambda$. If $\mu=w\cdot\lambda$ then there is functor 
$$
\gathered
I_w\,:\, \od_\GR(X)_{-\lambda}\ \longrightarrow\ \od_\GR(X)_{-\mu}\,,\ \ \text{such
that}
\\
\ \ (M\circ I_w)(\cf)\ = \ M(\cf)\,, \ \ \text{for}\ \ \cf \in
\od_\GR(X)_{-\lambda}\,.
\endgathered
\tag4.11
$$
The functors $I_w$ are called  intertwining functors and  were first introduced in 
\cite{BB2}. To give a formula for the functors $I_w$ we assume, for simplicity, that
$\lambda-\rho$ is integral. Let us set
$$
Y_w \ = \ \{\,(x,y)\in X\times X \mid \text{$y$ is in position $w$ with respect
to $x$}\,\}\,,
$$
and denote the projections to the first and the second factor by $p,q:Y_w\to X$,
respectively. The functor $I_w:\od(X)\to \od(X)$ is then given by:
$$
I_w(\cf) \ = \ Rq_*p^*(\cf)[\ell(w)]\,, \qquad \cf\in\od(X)\,,
$$
where $\ell(w)$ stands for the length of $w$.

\Lecture5{Matsuki correspondence for sheaves}			

In this lecture we explain geometric induction, the Beilinson-Bernstein localization,
and the Matsuki correspondence for sheaves.

\subhead{Geometric induction}\endsubhead
Let $A$ and $B$ be (linear) Lie groups such that $A\subset B$ and assume that
$B$   acts on  a semi-algebraic set $X$. We  construct a right adjoint 
$\Gamma_A^B$ to the forgetful functor $\operatorname{Forget}_{B}^A:\od_B(X)
\rightarrow \od_A(X)$. To this end, let us consider the diagram
$$\CD
X @<a<<B \times X @>q>> B/A \times X @>p>> X
\endCD
\tag5.1
$$
where $a(b,x)=b^{-1}x$, $q(b,x)= (bA,x)$, and $p(bA,x)=x$. The spaces in the diagram
have an action by $B\times A$ in such a way that the maps $a,q,p$ are 
$B\times A$-equivariant. This action of $B\times A$ on the spaces in (5.1) is given,
reading from left to right, by $(b,a)\cdot x = a\cdot x$\,,\
\ $(b,a)\cdot (b',x) = (bb'a^{-1},b\cdot x)$\,,\ \ $(b,a)\cdot (b'A,x) =
(bb'A,b\cdot x)$\,,\ \ $(b,a)\cdot x = b\cdot x$\,. To give a formula for the functor
$\Gamma_A^B$, let us pick $\cf\in \od_A(X)$. As the $B$-action is trivial on $X$, we can
view $\cf\in \od_{B\times A}(X)$. Then, by property (3) of the characterization of
the equivariant derived category, there is a unique $\tilde\cf\in \od_B( B/A \times
X)$ such that $q^*\tilde{\F} = a^*\F$. We then set
$$
\Gamma^B_A\F = Rp_* \tilde\F\,, \qquad \text{where $\tilde\F$ is the unique sheaf
such that 
$q^*\tilde{\F} = a^*\F$.}
\tag5.2
$$
Intuitively, the operation $\Gamma^B_A$ amounts to averaging $\cf$ over
''$B/A$-orbits''.

\subsubhead{Parabolic induction}\endsubsubhead
We pause briefly to explain how to phrase  parabolic induction in terms of the
geometric induction functors. For simplicity, we do it for the
trivial infinitesimal character only. Let
$P_\RR \subset \GR$ be a parabolic subgroup with Levi decomposition  $P_\RR = L_\RR
N_\RR$ with $P=LN$ the corresponding complexified Levi decomposition. We denote by
$X_L$  the flag manifold of the group $L$. Associated to $\F\in\od_{L_\RR}(X_{L})$
we have a (virtual) representation $M(\cf)$ of $L_\RR$. In the fibration $X
\rightarrow G/P$ the fiber over the point $eP$ can be identified with $X_L$ and we
denote by $i:X_L \hookrightarrow X$ the inclusion. Then we have the following formula
for parabolic induction:
$$
  \Ind_{P_\RR}^{G_\RR}(M(\F)) = M(\Gamma_{L_\RR}^\GR i_* \F) .
\tag5.3
$$

\subhead Beilinson-Bernstein localization\endsubhead

In lecture 4 we explained how to associate a
$\GR$-representation $M(\cf)$ to an element $\cf\in\od_\GR(X)_{-\lambda}$. We will now
explain how to construct the Harish-Chandra module associated to $M(\cf)$. For this
we fix a maximal compact subgroup $\KR$ of $\GR$ and denote by $K$ the
complexification of
$\KR$. The answer is provided by the following commutative diagram:
$$
\CD
\{\text{$\GR$-representations}\}_{\chi_\lambda} @>>>
\{\text{H-C-modules}\}_{\chi_\lambda}
\\
@AMAA @AA{\alpha}A
\\
\od_\GR(X)_{-\lambda} @>\Gamma>> \od_K(X)_{-\lambda}\,.
\endCD
\tag5.4
$$ 
The arrow on the top row associates to a representation its Harish-Chandra module. The
arrow $\alpha$, due to Beilinson and Bernstein \cite{BB1}, amounts to taking the
cohomology of the $\Cal D$-module that is associated  to the element in
$\od_K(X)_{-\lambda}$ by the Riemann-Hilbert correspondence\footnote{For a
treatment of $\Cal D$-modules and the Riemann-Hilbert correspondence, see
\cite{Bo}}. One can write the functor
$\alpha$, in analogy with
$M$, as
$$
\alpha\ :\  \F \longmapsto \sum(-1)^k \Ext^k(\Dual\cf,\O^{\text{alg}}(\lambda))\,, 
\tag5.5
$$
where $\O^{\text{alg}}(\lambda)$ is  the sheaf of twisted algebraic functions on
$X$. It is a subsheaf of $\O(\lambda)$. Note that a similar analogue of the
functor $m$ does not make sense, i.e., it does not produce a Harish-Chandra module.
Finally, as to the functor $\Gamma$,
$$
\Gamma\,:\,\od_\GR(X)_{-\lambda}\  @>{\ \sim \ }>> \  \od_K(X)_{-\lambda}\,,\qquad
\Gamma = \Gamma_{K_\RR}^K \circ \operatorname{Forget}_{K_\RR}^\GR
\tag5.6a
$$
is an equivalence of categories \cite{MUV}. The functor $\Gamma$  has the following
property which justifies calling it the Matsuki correspondence for sheaves:
$$
\Gamma(Rj_!\CC_{\O})\  = \ Rj'_*\CC_{\O'} [-2 \operatorname{codim}_\CC \O']
\tag5.6b
$$ 
if the $\GR$-orbit $\Cal O$ corresponds to the $K$-orbit $\Cal O'$ under the Matsuki
correspondence, which we recall below in (5.7). Here $j:\Cal O \hookrightarrow X$ and
$j':\Cal O'
\hookrightarrow X$ denote the inclusions of the orbits $\Cal O$ and $\Cal O'$ to the
flag manifold $X$. We will explain in some detail below the geometric idea behind the
proof of (5.6).

\subheading{Matsuki Correspondence for sheaves} 
The main ingredient of the proof of the Matsuki correspondence for sheaves 
is a Morse theoretic interpretation and refinement of the original result of Matsuki.
Let us recall Matsuki's statement: there is  a bijection
$$
  \GR \backslash X \ \longleftrightarrow K\backslash X
\tag5.7
$$
between $\GR$-orbits on $X$ and $K$-orbits on $X$ such that  a $\GR$-orbit $\O'$
corresponds to a $K$-orbit $\O$ if and only if  $\O' \cap \O$ is non-empty and
compact. Furthermore, $\GR$-equivariant local systems on $\O'$ correspond
bijectively to $K$-equivariant local systems on $\O$. The fundamental idea, due
to Uzawa, is that there exists a Bott-Morse function $f$ on $X$ whose stable
manifolds, with
respect to a particular metric,  are  the $K$-orbits and whose unstable manifolds are
the $\GR$-orbits.   To construct the metric and the function, we write $\gr=\kr\oplus
\fp_\Bbb R$ for the Cartan decomposition, and let $U_\RR$ denote the compact form
corresponding to  $\frak u_\RR = \frak k_\RR \oplus i \frak p_\RR$.  We choose an
$H$-integral, regular, dominant $\lambda\in\fh^*$, and denote the corresponding
highest weight representation of $G$ by $V_\lambda$. Then the line bundle $\bold
L_\lambda$ gives us an embedding $X\hookrightarrow \Bbb P(V_\lambda)$. We fix a
$\UR$-invariant Hermitian scalar product on $V_\lambda$. The real part of the scalar
product induces a $\UR$-invariant Riemannian metric, the Fubini-Study metric on the
projective space $\Bbb P(V_\lambda)$ and hence on $X$. This is the metric we use.

To construct the desired function $f$, we note that the element $\lambda\in\fh^*$ gives
us an embedding
$$
X \ @>{\ \sim \ }>> \Omega_\lambda \subset i\frak u_\Bbb R^*\,, \qquad
\Omega_\lambda\ \text{a $\UR$-orbit}
\tag5.8
$$ 
as follows. Given $x\in X$, there is a unique Cartan $\TR$ of $\UR$ which fixes $x$.
As explained in lecture 1, this gives us a map $\tau_x :\ft_\Bbb R \to \fh$ which lifts
to a map
$\TR\to H$ and thus $\lambda\in\fh^*$ gives rise to  an element 
$\lambda_x\in i\ft^*_\Bbb R$. Via the direct sum decomposition
$\frak u_\Bbb R= \ft_\Bbb R\oplus [\ft_\Bbb R,\frak u_\Bbb R]$, we can view
$\lambda_x\in i\frak u_\Bbb R^*$ and the association $x\mapsto \lambda_x$ gives a
map from $X$ to $i\frak u_\Bbb R^*$. As $ i\frak u_\Bbb R^*= i\frak k_\Bbb R^*\oplus
\frak p_\Bbb R^*$, we get a map
$$
m\,:\, X \longrightarrow i\frak u_\Bbb R^*\longrightarrow\frak p_\Bbb R^*\,.
\tag5.9
$$
The Killing form induces a metric on $\frak p_\Bbb R$ and hence on its dual $\frak
p_\Bbb R^*$. We define
$$
f\,:\, X \longrightarrow \Bbb R\,, \qquad \text{by} \ \ f(x)\ = \ \|m(x)\|^2\,.
$$
Here is the refined version of the Matsuki correspondence:
$$
\gathered
\text{The function $f$ is a Bott-Morse function. Its gradient flow $\nabla f$,} 
\\
\text{with respect to the metric on $X$ described above, has the}
\\
\text{$K$-orbits as stable manifolds and $\GR$-orbits as unstable manifolds.}
\\
\text{The critical set consists of a finite set of $\KR$-orbits.}
\endgathered
\tag5.10
$$
The statement (5.10) is illustrated below in figure 1 for the case $\GR=SL(2,\Bbb R)$.
\midinsert
\centerline{\BoxedEPSF{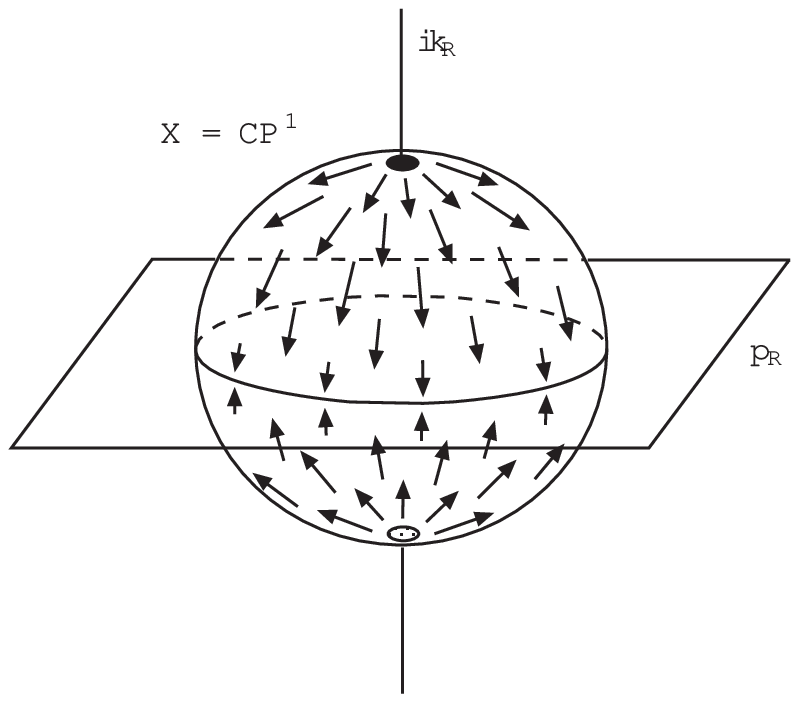}}
\botcaption {Figure~1} \endcaption
\endinsert

In proving the Matsuki correspondence (5.6) for sheaves it is crucial to know that if
a
$K$-orbit $\Cal O$ and a $\GR$-orbit $\Cal O'$ are related by the Matsuki
correspondence, then both $\Cal O$ and $\Cal O'$ can be retracted to $\Cal O\cap\Cal
O'$. This is the content of statement (5.10).

\remark{Remark} In  \cite{N}, Ness defined a ''moment map" $\mu: \Bbb
P(V_\lambda) \to i\frak u_\Bbb R^*$ for the action of the compact group $\UR$
on $\Bbb P(V_\lambda)$. The map (5.8) is the restriction of this moment map to $X$
via the embedding $X \hookrightarrow \Bbb P(V_\lambda)$. The construction of Ness
readily extends to actions of  semisimple groups, see, for example, \cite{MUV}. The
moment map for the $\GR$-action on  $\Bbb P(V_\lambda)$ composed with the the embedding
$X \hookrightarrow \Bbb P(V_\lambda)$ is exactly the map $m$ in (5.9). 
\endremark

\subhead Representations and perverse sheaves \endsubhead

In these lectures we emphasize the role of $\od_\GR(X)_{-\lambda}$ over
$\od_K(X)_{-\lambda}$. However, for certain things it is preferable to work with the
category  $\od_K(X)_{-\lambda}$ instead, as we will explain below. Let
$\lambda\in\fh^*$ be dominant and let
$\operatorname{P}_K(X)_{-\lambda}$ denote the  subcategory of $\od_K(X)_{-\lambda}$
of perverse sheaves \cite{BBD}. If $\cf\in\operatorname{P}_K(X)_{-\lambda}$ then, as
was shown in \cite{BB1} in the language of $\Cal D$-modules,
$$
 \Ext^k(\Dual\F, \O(\lambda)) \ = \ 0\,, \qquad \text{if} \ \ k\neq 0\,.
$$
Thus, the  functor $\alpha$ of (5.5), restricted to
$\operatorname{P}_K(X)_{-\lambda}$, gives a functor
$$
  \alpha\ : \ \operatorname{P}_K(X)_{-\lambda}\ \longrightarrow \ 
\{\text{Harish-Chandra modules}\}_{\chi_\lambda}\,,
\tag5.11
$$
which is an equivalence of categories if $\lambda$ is regular. If $\lambda$ is not
regular then there is a kernel  which can be described explicitly \cite{BB1}. Under
(5.11) the irreducible Harish-Chandra modules correspond to intersection homology
sheaves of $K$-equivariant irreducible local systems on $K$-orbits. 

Note that the category $ \operatorname{P}_K(X)_{-\lambda}$, as well as intersection
homology sheaves, are characterized by conditions on objects of
$\od_K(X)_{-\lambda}$ which are local on $X$.  This does not appear to be possible on
the $\GR$-side. The nice subcategory analogous to  $
\operatorname{P}_K(X)_{-\lambda}$ exists on the $\GR$-side for formal reasons. By
the commutativity of (5.4) we can simply take it to be 
$\Gamma^{-1}(\operatorname{P}_K(X)_{-\lambda})$. The functor
$\gamma=\Gamma^{-1}:\od_K(X)_{-\lambda}\to \od_\GR(X)_{-\lambda}$ is given
analogously to
$\Gamma$, by switching the roles of $K$ and $\GR$, and by replacing all the *'s in
the construction by !'s. From this description one can see  that
there can not be a characterization of $\Gamma^{-1}(\operatorname{P}_K(X)_{-\lambda})$ 
as a subcategory  of $\od_\GR(X)_{-\lambda}$  using
conditions  which are local on $X$. 

\subsubhead{Cohomological induction}\endsubsubhead
Finally, let us describe cohomological induction in geometric terms. For simplicity,
we do so only in the case of trivial infinitesimal character. Let $P\subset G$ be a
parabolic subgroup such that its Levi $L$ is $\theta$-stable and defined over $\RR$,
i.e., $L_\RR \subset G_\RR$ and $L_\RR \cap K_\RR \subset L_\RR$ is maximal compact.
Let us consider the fibration $X\to G/P$. Its fiber over $eP$ can be identified with
the flag manifold of $X_L$ of $L$, and we denote the inclusion of that fiber in $X$
by $i:X_L \to X$. Then
$$
  \Ind^{(\fg,K)}_{(\frak l, L\cap K)} \alpha(\F) = \alpha(\Gamma^K_{L\cap K} i_*
\F)\ ,
$$
where $\Ind^{(\fg,K)}_{(\frak l, L\cap K)}$ stands for the cohomological induction (in
the sense of \cite{EW}) from $(\frak l, L\cap K)$-modules to $(\fg,K)$-modules. For
this fact see, for example, \cite{MP} and \cite{S2}.

\Lecture6{Characteristic cycles}			

Let $X$ be a real algebraic manifold. 
For simplicity, we assume that $X$ is oriented.
As usual, let $\od(X)$ be the bounded derived category of constructible sheaves
on $X$. 
A simple invariant that one can associate to an object $\F \in \od(X)$ is
its {\it Euler characteristic}, $\chi(X,\F)$, defined by
$$
  \chi(X,\F) = \sum_{k} (-1)^k \dim_{\CC} \oh^k(X,\F)\, .
\tag6.1
$$
The Euler characteristic is additive in distinguished triangles (recall that
triangles arise from exact sequences of complexes):  if 
$$
\F' \rightarrow \F \rightarrow \F''  \rightarrow \F'[1]
\tag6.2
$$ 
is a distinguished triangle then 
$$
\chi(X,\F) = \chi(X,\F') + \chi(X,\F'')\, .
\tag6.3
$$ 
As a  particular special case we get $\chi(X,\cf[1])=-\chi(X,\cf)$.
From (6.2) we also conclude that $\chi(X,\ \, ): \od(X)
\rightarrow \Bbb Z$ descends to the
$K$-group $K(\od(X))$. Recall that the $K$-group is the free abelian group
generated by all the $\cf\in\od(X)$ subject to the relation $\cf=\cf'+\cf''$ for
all distinguished triangles (6.2). The local Euler characteristic
$$
  \chi(\F)\ :\ X \ \longrightarrow\ \ZZ
\tag6.4
$$
is defined  by 
$$
  \chi(\F)_x \ = \ \sum_{k} (-1)^k \dim_{\CC} H^k(\F)_x\,.
\tag6.5
$$
The local Euler characteristic $\chi(\cf)$ is a constructible function: there is a
triangulation $\T$ of $X$ such that for any $\sigma\in\T$ the restriction
$\chi(\cf)|\sigma$ to the open simplex $\sigma$ is constant. 

\remark{Remark} It is not very difficult to verify that the homomorphism (both sides
are  abelian groups) 
$$ 
\matrix
  \operatorname K(\od(X)) & \longrightarrow & \{\text{constructible functions}
\ X \rightarrow \ZZ \, \} \\
 & & \\
\F & \longmapsto & \chi(\F)
\endmatrix
\tag6.6
$$
is an isomorphism.
\endremark

To give a more geometric description of the category $K(\od(X))$, we recall
another point of view to the Euler characteristic.  Assume now that  $X$ is
compact. Then, by a classical theorem of Hopf, the Euler characteristic $\chi(X) =
[X].[X]$, the self intersection product of the zero section $[X]$ in
$T^*X \simeq TX$. To be able to take the self intersection product, one perturbes
one copy of the zero section so that we stay in the same homology class. For
example, we can replace one of the copies of $X$ by a generic section of the
tangent bundle. 

To generalize the theorem of Hopf to arbitrary $\cf\in\od(X)$, Kashiwara, in
\cite{K1}, introduced the  notion of a {\it characteristic cycle\/} $\cc(\F)$. The
characteristic cycle  $\cc(\F)$ is a  semi-algebraic Lagrangian cycle (not necessarily
with compact support) on the  cotangent bundle $T^*X$. The definition of
$\cc(\F)$ is Morse-theoretic. Heuristically, $\cc(\F)$ encodes the infinitesimal
change of the local Euler characteristic $\chi(\cf)_x$ to  various
co-directions in $X$. We have the following generalization of Hopf's theorem
\cite{K1}:
$$
\chi(X,\cf) \ = \ \occ(\cf).[X]\,.
\tag6.7
$$
As the definition of $\occ$ is quite technical, we will omit it, and refer to
\cite{K1} and \cite{KS, chapter IX} for details. Below we will give an axiomatic
characterization of $\occ$ following \cite{SV1}.

\subsubhead{Semi-algebraic chains and cycles}\endsubsubhead
Let $M$ be a real algebraic manifold (In our situation, $M= T^*X$).
A {\it semi-algebraic $p$-chain\/}, $C$, on $M$ is  a finite 
integer  linear combination
$$
C\ = \ \sum n_\alpha[S_\alpha]\,.
\tag6.8
$$  
Here the $S_\alpha$ are oriented $p$-dimensional semi-algebraic submanifolds of $M$
and the symbols $[S_\alpha]$ are subject to  the following relations:

\roster
\item"(i)" $[S_1\cup S_2] = [S_1] + [S_2]$ if $S_1,S_2$ are disjoint;
\item"(ii)" $[S^-] = -[S]$, where $S^-$ is the manifold $S$ with the 
opposite orientation;
\item"(iii)" $[S] = [S']$ if $S'\subset S$ is an open and dense subset of $S$
with the orientation induced from $S$.
\endroster

From the relations above we conclude that we can write any chain $C$ in (6.8) in
such a way that the $S_\alpha$ are disjoint. Once we have done so, the support of
 $C$ is defined as $|C|  = \overline{\bigcup_{\alpha} S_\alpha}$, the closure of the
union of the $S_\alpha$. As semialgebraic sets can be triangulated, we can define
the boundary operator $\partial$ from $p$-chains to  $(p-1)$-chains in
the usual way. If $C$ is a $p$-chain and $\partial C = 0$ then we call
$C$ a $p$-cycle.

\subsubhead{Lagrangian cycles}\endsubsubhead 
Let $M=T^*X$, where $X$ is a real algebraic  manifold of dimension $n$. The
manifold $\ct$ has a canonical symplectic structure. We call a semi-algebraic
subset $Z$ of $\ct$ Lagrangian if $Z$ has an open dense subset $U$ consisting of
smooth points such that $U$ is a Lagrangian submanifold of $\ct$. We call a cycle
on $\ct$ Lagrangian if its support is. The group of positive reals $\Bbb{R}^+$ acts
by scaling on $\ct$.  We denote by $\Cal{L}^+(X)$ the group of semi-algebraic,
$\Bbb{R}^+$-invariant Lagrangian  cycles on $T^*X$. Each $C\in\Cal{L}^+(X)$ is an
$n$-cycle on $\ct$ and, as is not very hard to show,
$$
|C| \ \subset \ \bigcup T^*_{S_i}X\,,  \qquad  S_1, \dots S_k\subset X\ \ \ 
\text{submanifolds}\,.
\tag6.9
$$

\subhead Characteristic Cycles \endsubhead

We  will now give the axiomatic description of the characteristic cycle
construction. For that we assume, for simplicity, that the real algebraic manifold
$X$ is orientable, and we fix an orientation of $X$. The characteristic cycle
construction $\occ$ is a map
$$
\cc\,:\, \od(X) \longrightarrow \Cal{L}^+(X)
\tag6.10
$$
satisfying the following  properties: 

\roster 
\item"(a)" The definition is local, {\it i.e.,} the following diagram commutes:
$$
\CD
\od(X) @>{\ \ \ \occ\ \ \ }>> \Cal L^+(X)  
\\
@V{j^*}VV              @VV{j^*}V   
\\
\od(U)  @>{\ \ \ \occ\ \ \ }>>  \Cal L^+(U)\,;
\endCD
$$
here $U$ is an open subset of $X$ and the $j^*$ on the right denotes the
restriction\footnote{The restriction can be performed because our cycles do not
necessarily have compact support.} of  cycles from
$\ct$ to
$T^*U$. 
\item"(b)" $\cc(\Bbb{C}_X) = [X]$; the symbol $[X]$ makes sense because we have
fixed an orientation of $X$. 
\item"(c)" $\cc$ is additive in exact sequences,  i.e., if 
$\F' \rightarrow \F \rightarrow \F''  \rightarrow \F'[1]$
is a distinguished triangle then 
$$
\cc(\F)\ = \ \cc(\F') \ +\  \cc(\F'')\,.
$$
\item"(d)" 
If $j: U \rightarrow X$ is an open embedding, and $g$ is any semi-algebraic
defining equation of $\partial U$  -- which is at least $C^1$ and which we assume to be
positive on $U$ --   then the following diagram commutes:
$$
\CD
\od( U) @>{\ \ \ \occ\ \ \ }>> \Cal L^+(U)  \ \ 
\\ @V{Rj_*}VV              @VV{j_*}V   
\\
\od( X)  @>{\ \ \ \occ\ \ \ }>>  \Cal L^+(X) \ ;
\endCD
$$
here $j_*: \Cal L^+(U) \to  \Cal L^+(X)$ is the following limit operation
$$
j_*(C) = \lim_{s \to 0^+} (C + s\, d\,\text{log}\,g)\,, \ \ \ C \in \Cal 
L^+(U)\,.
$$
\endroster

In property (d) the notation $C + s\, d\,\text{log}\,g$ stands for the cycle which
is obtained by applying the automorphism $(x,\xi)\mapsto
(x,\xi+s\frac {dg_x}{g(x)})$ of $T^*U$ to $C$. The limit operation can be
interpreted as
$$
 \lim_{s \to 0^+} (C + s\, d\,\text{log}\,g) \ = \ -\partial \tilde C\,,
$$
where we view
$$
\tilde C \ = \ \{C + s\, d\,\text{log}\,g\mid s>0\}
$$
as an $(n+1)$-chain in $\Bbb R\times\ct$. The property (4) is proved as a theorem
in \cite{SV1}, where one can also find a more detailed discussion of the notion of
limit and precise conventions for orientations. 

\definition{Example 1}
Let $j: (0,1) \hookrightarrow \RR$ be the inclusion map of the open interval. As a
defining equation of $\partial (0,1)$ we choose  $g(x) = x(1-x)$. We will apply
(4) with $U=(0,1)$, $X=\Bbb R$, and $\cf= \Bbb C_{(0,1)}$. Because, by (2), 
$\occ(\cf)=[(0,1)]$, we get
$$
\split
\cc(Rj_*\CC_{(0,1)}) &= \lim_{s\rightarrow
0^+}\left[\left\{s\frac{dg_x}{g(x)} \mid 0<x<1\right\}\right]
\\ &= \lim_{s\rightarrow
0^+}\left[\left\{s\frac{dx}{x}-s\frac{dx}{1-x} \mid 0<x<1\right\}\right].
\endsplit
$$
The result of this calculation is illustrated below in figure 2.
\midinsert
\centerline{\BoxedEPSF{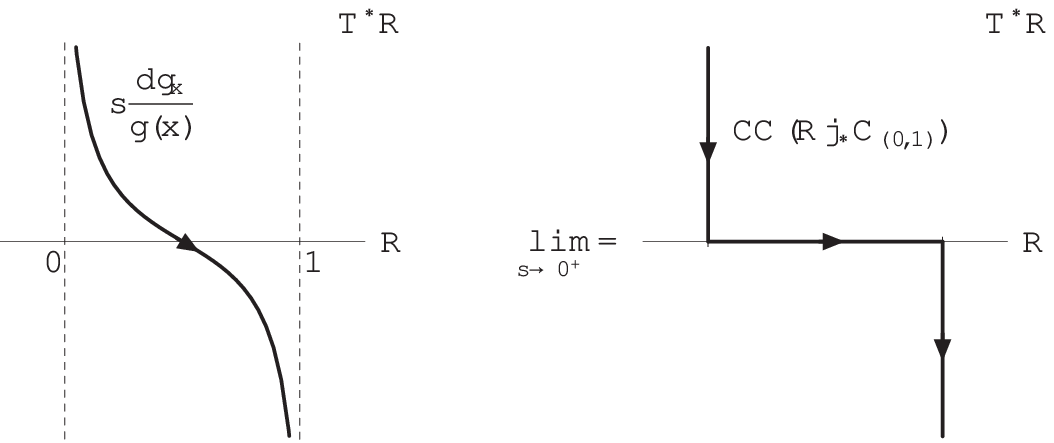}}
\botcaption {Figure~2} \endcaption
\endinsert
\enddefinition

Proceeding in the same manner as in the previous example, we see that
$$
\occ(\Bbb C_Y) \ = \ [T^*_YX], \qquad Y \ \text{a closed submanifold of $X$}\,,
$$
for a particular orientation (see the orientation conventions of \cite{SV1}) of
$T^*_YX$.

\subhead{Further properties of $\occ$}\endsubhead

Kashiwara's theorem (6.7) extends to the relative situation. Let $f:X \rightarrow Y$
be a proper  real algebraic map between (oriented) real algebraic manifolds.
Then 
$$
\cc(Rf_* \F)\  = \ f_*(\cc(\F)).
\tag6.11
$$
To define the map 
$$
  f_*\,:\, \Cal{L}^+(X) \ \longrightarrow \ \Cal{L}^+(Y)
\tag6.12
$$
on the right hand side of the equation, let us
consider the commutative diagram 
$$
\CD 
T^*X @<{df}<< X \times_Y T^*Y @>\tau>> T^*Y \\
@V{\pi_X}VV  @VVV  @VV{\pi_Y}V \\
X @<=<< X @>f>> Y
\endCD 
\tag6.13
$$
The assumption that $f$ is proper implies that $\tau$ is proper. The map (6.12) is
defined by the formula
$$
f_*(C) \ = \ \tau_*((df)^*C)\,.
\tag6.14
$$
Here $(df)^*$ is a Gysin map. If the cycle $C$ happens to be transverse to the map
$df$ then $(df)^*(C)= (df)^{-1}(C)$. If this is not the case, then we embed $C$
into a family of cycles whose generic members $C_s$ are transverse to $df$. In this
way $(df)^*(C)$ is expressed as the limit $(df)^*(C) = \lim_{s\to 0^+}
(df)^{-1}(C_s)$. For a proof of (6.11), see \cite{KSa} and for its interpretation
in the present language see \cite{SV1}. 

Finally, let us describe the effect of the operation $f^*$ on characteristic cycles
for $f:X\to Y$ a submersion. By property (a) it suffices to do so when $X=Z\times
Y$ and $f:Z\times Y\to Y$  is the projection. Then
$$
\occ(f^*\cg) \ = \ \occ (\Bbb C_Z \boxtimes \cg) \ = \ [Z]\times\occ(\cg)\,.
\tag6.15
$$
This formula follows from the properties (a-d). On the contrary, there does not
appear to be an easy formal way to deduce (6.11) from properties (a-d). Finally,
from property (c) we conclude that
$$
\occ(\cf) \ = \ \tsize\sum (-1)^k \occ(\oh^k(\cf))
\tag 6.16
$$
for $\cf\in\od(X)$.

\definition{Example} Let us deduce (6.7) from (6.11). To this
end, let us assume that $X$ is a compact (oriented!) real algebraic manifold,
$Y=\{\pt\}$, and $f:X\to\{\pt\}$. applying (6.16) to $Rf_*\cf$, gives     
$\cc(Rf_* \F) = \chi(\F)$. By formula (6.11) we get:
$$
\chi(\F) \ = \ \cc(Rf_* \F) \ = \ \tau_*(df)^*(\cc(\F)) \ = \ \occ(\cf).[X]\,;\,
$$
with an appropriate interpretation of the signs in the intersection product.  
\enddefinition

\Lecture7{The character formula}			

In this lecture we give an integral formula for the Lie algebra character of the
representation $M(\cf)$ associated to a $\cf\in\od_\GR(X)_{-\lambda}$ in terms of the
characteristic cycle $\occ(\cf)$. We begin by briefly recalling the notion of the
character on the Lie algebra.

\subheading{The Lie algebra character}

Let $\pi$ be a representation, in our previous sense, of the semisimple Lie group
$\GR$ with infinitesimal character $\chi_\lambda$. To $\pi$, following Harish-Chandra,
we can associate its character $\Theta_\pi$. The character  $\Theta_\pi$ is an
invariant eigendistribution on $\GR$, i.e., it is conjugation invariant and the center
$\zg$ of the universal enveloping algebra $\ug$ acts on $\Theta_\pi$ via the character
$\chi_\lambda$. Via the exponential map we can, at least in the neighborhood of the
origin, pull back the distribution $\Theta_\pi$ to the Lie algebra $\gr$. We define
the Lie algebra character by the formula
$$
\tp =
\sqrt{\operatorname{det}(\operatorname{exp}_*)}\,\operatorname{exp}^*\TP\,.
\tag7.1
$$
We have inserted the factor $\sqrt{\operatorname{det}(\operatorname{exp}_*)}$ so that
$\tp$ is an invariant eigendistribution on $\gr$, i.e., conjugation invariant and the
constant coefficient differential operators $S(\fg)^G\cong\zg$ act on $\tp$ via the
character $\chi_\lambda$. Harish-Chandra's regularity theorem implies that $\tp$ can
be extended from the neighborhood of the origin uniquely to all of $\gr$ and that $\tp$
is a locally $L^1$-function which is real analytic on the set of regular semisimple
elements in
$\gr$. 

\subheading{Rossmann's formula}

Let us now assume that the group $\GR$ has a compact Cartan. As
was explained in lecture 4, the discrete series representations are parametrized by
regular, antidominant $\lambda\in \Lambda+\rho$ and a choice of an open $\GR$-orbit on
$X$. To this data we attach a discrete series representation $\pi=\pi(S,\lambda)$.
Recall (4.5): as virtual representations, $\pi(S,\lambda)= (-1)^sM(Rj_*(\Bbb
C_S))$, where $j:S\hookrightarrow X$ denotes the inclusion.  Associated to the data
$(S,\lambda)$ we construct a coadjoint $\GR$-orbit
$\Omega_\lambda(S)$ in $i\fg_\Bbb R^*$ exactly the same way as in lecture 1 where we
had assumed that $\GR$ is compact. Given $x\in S$, there is a unique compact Cartan
$\TR$ that fixes $x$. This gives us a map $\tau_x:\frak t_\Bbb R \to \fh$ which, in
turn, allows us to view
$\lambda$ as an element $\lambda_x\in i\frak t_\Bbb R^*$. Finally, the direct sum
decomposition $\gr = \frak t_\Bbb R \oplus [\frak t_\Bbb R,\gr]$, allows us to view
$\lambda_x\in i\fg_\Bbb R^*$. The association $x\mapsto \lambda_x$ gives a
$\GR$-equivariant identification of $S$ with a $\GR$-orbit $\Omega_\lambda(S)\subset
i\fg_\Bbb R^*$. The orbit $\Omega_\lambda(S)$ has a canonical symplectic form
$\sigma_\lambda$.
Let $\theta_\pi$ denote the Lie algebra character of the discrete series
representation $\pi(S,\lambda)$. Rossmann, in \cite{R1}, proves the following result:
$$
\int_{\gg_\RR}
\theta_\pi\, \phi\, dx \ = \ \frac 1 {(2\pi i)^n n!}
\int_{\Omega_\lambda(S)}  \hat \phi \ \sigma_\lambda^n
\,,
\tag7.2
$$
in complete analogy with (1.6), where $\phi$ is a tempered test function and the Fourier
transform is performed without the $i$ as in (1.5). As in lecture 1, this result can be
phrased as
$$
\hat \theta_\pi\ = \ \text{the coadjoint orbit
$\Omega_\lambda(S)$ with measure $\frac{\sigma_\lambda}{(2\pi i)^n n!}$}\,.
\tag7.3
$$
From this one can obtain  an analogous formula for all tempered
representations. 

Formula (7.2,3), as stated, can not be generalized for non-tempered representations:
the Fourier transform $\hat\theta_\pi$ no longer makes sense when $\tp$ is not
tempered. However, in \cite{R2} Rossmann proposed a way to generalize the integral
formula (7.2). More specifically, he showed how to write down the  invariant
eigendistributions on $\gr$  as integrals resembling (7.2). We begin by
constructing the twisted moment map
$$
\mu_\lambda\,:\,\ct \ \longrightarrow \fg^*\,.
\tag7.4
$$
We consider the compact form $\UR$ corresponding to the Lie algebra $\frak u_\Bbb R =
\kr\oplus i\fp_\Bbb R$. Using the same construction as in (5.8), we obtain a 
$\UR$-equivariant real algebraic map
$$
m_\lambda\,:\, X \  \longrightarrow i\frak u_\Bbb R^* \subset \fg^*\,.
\tag7.5
$$
Recall the moment map 
$$
\mu:\ct \to \fg^*\,,
\tag7.6
$$
which is given, on the level of fibers 
$T^*_xX\cong (\fg/\fb_x)^*$, by the canonical inclusion $(\fg/\fb_x)^* \to \fg^*$. The
moment map $\mu$ is $G$-equivariant and complex algebraic. Its image is the nilpotent
cone, once we make  the identification $\fg\cong \fg^*$ by means of the Killing form.
The twisted moment map $\mu_\lambda$ is given by the following formula:
$$
\mu_\lambda \ = \ \mu  \ + \ m\circ \pi\,,
\tag7.7
$$
where $\pi:\ct \to X$ denotes the projection. The twisted moment map, for regular
$\lambda$, provides an isomorphism 
$$
\mu_\lambda\,:\,\ct @>{\ \ \sim \ \ }>> \ \Omega_\lambda,
\tag7.8
$$
where $\Omega_\lambda\subset \fg^*$ is a $G$-orbit. Let us denote the canonical
symplectic form on $\Omega_\lambda$ by $\sigma_\lambda$ and let us write
$$
T^*_\GR X \ = \ \bigcup_{\text{$S$ a $\GR$-orbit}} T^*_SX\,.
\tag7.9
$$
Then we see that
$$
\ohf {2n}(T^*_\GR X, \Bbb C)\ = \ \{\text{Lagrangian $\Bbb C$-cycles on $\ct$ supported
on $T^*_\GR X$}\}\,,
\tag7.10
$$
where $n=\dim_\Bbb C X$, and  the symbol $\ohf {2k}$ stands for homology with closed,
i.e., possibly non-compact, supports. 
We define the Fourier transform $\hat\phi$ of a test
function $\phi$ in $C^\infty_c(\gr)$ without choosing a square root of $-1$\,,
as a holomorphic function on $\fg^*$\,:
$$
\hat \phi (\zeta) \ = \ \int_\gr e^{\zeta(x)}\phi(x) dx \qquad
(\zeta\in\fg^*)\,.
\tag7.11
$$
Let us
assume that $\lambda$ is regular. Rossmann \cite{R2} shows that
$$
\gathered
\ohf {2n}(T^*_\GR X, \Bbb C) @>{\ \ \sim\ \ }>>
\left\{\foldedtext\foldedwidth{50\jot}{invariant eigendistributions on $\gr$ with
infinitesimal character $\chi_\lambda$}\right\}
\\
C\ \ \ \mapsto\ \ \ \ \left\{\phi \mapsto  \frac 1 {(2\pi i)^n n!}
\int_{\mu_\lambda(C)}  \hat \phi \ \sigma_\lambda^n \right\}\,.
\endgathered
\tag7.12
$$
Here $\phi\in C^\infty_c(\gr)$ is a test function. The integral
$\int_{\mu_\lambda(C)}  \hat \phi \ \sigma_\lambda^n$ converges because $\hat\phi$
decays rapidly in the imaginary directions and the cycle $\mu_\lambda(C)$ has bounded
real parts: $\mu_\lambda(C)$ differs from $\mu(C)$ by a compact set and $\mu(T^*_\GR
X)\subset
i\fg_\Bbb R^*$. Note that the form $\hat\phi\sigma_\lambda^n$ is holomorphic,
of top degree in $\Omega_\lambda$. Hence, the cycle $\mu_\lambda(C)$   can be replaced
by a homologous cycle without changing the value of the integral
$\int_{\mu_\lambda(C)}  \hat \phi \ \sigma_\lambda^n$\,, provided that the chain
giving rise to the homology has bounded real parts.

\subheading {The character formula} 
Fix $\F\in \od_\GR(X)_{-\lambda}$ and let $\theta(\F)$ denote the Lie algebra character
 of the representation $M(\cf)$. We continue to assume that $\lambda$ is regular. Then,
as is shown in \cite{SV2}:
\proclaim{Theorem} The character $\theta(\cf)$ is given by taking $C=\occ(\cf)$ in
(7.12), i.e., 
$$
  \int_{\fg_\RR} \theta(\F) \phi dx \ = \ \frac{1}{(2\pi i)^n n!}
  \int_{\mu_\lambda(\cc(\F))} \hat{\phi}\sigma_\lambda^n\,,
\tag7.13
$$
for $\phi\in C^\infty_c(\gr)$.
\endproclaim
\remark {Remark} The fact that $\occ(\cf)$ is supported on $T^*_\GR X$ follows from
the $\GR$-equi\-variance of $\cf$. Note that $\occ(\cf)\in \ohf {2n}(T^*_\GR X, \Bbb
Z)$, i.e., the characters are given by integral cycles in (7.12). 
\endremark

To extend the validity of the theorem to arbitrary $\lambda$, we use the following
result:
\proclaim{Lemma}
$\mu_\lambda^* \sigma_\lambda = - \sigma + \pi^*\tau_\lambda$, where
$\sigma$ denotes the canonical symplectic form on $T^*X$ and $\tau_\lambda$ 
is a $2$-form on $X$ defined by $\tau_\lambda(u_x,v_x) = \lambda[u,v]$, $x\in X$, and
 $u_x,v_x$ are tangent vectors  given by $u,v \in \frak u_\RR$.
\endproclaim

The following formula is valid for all $\lambda$:
$$
  \int_{\fg_\RR} \theta(\F) \phi dx \ = \ \frac{1}{(2\pi i)^n n!}
  \int_{\cc(\F)} \mu_\lambda^* \hat{\phi} ( - \sigma + \pi^*\tau_\lambda)^n\,.
\tag7.14
$$
 
As an example, let us consider the case $\GR=SL(2,\Bbb R)$. Then $X=\Bbb C\Bbb P^1$,
and there are three $\GR$-orbits, the upper and lower hemispheres $D_+$, $D_-$, and
the equator $\Bbb R \Bbb P^1\cong S^1$. Let us consider a discrete series
representation associated to $D_-$ and a negative $\lambda \in\Bbb Z\cong\Lambda$. We
take $\cf = Rj_*\CC_{D_-} = \CC_{\bar{D}_-}$. To calculate $\cc(Rj_*\CC_{D_-})$, we
view $D_-$ as a unit disc in $\Bbb C$ and, by the defining property (d) of the
characteristic cycle in lecture 6, we should choose a defining equation $f$ for the
boundary of $D_-$. Once the equation is chosen, we get:
$$
\aligned
\cc(Rj_*\CC_{D_-}) &= \lim_{s \rightarrow 0^+}\left(\occ(\CC_{D_-}) +s
\frac{df}{f}\right)\\ &=  \lim_{s \rightarrow 0^+}\left\{s \frac{df_x}{f(x)} 
\mid x \in
D_+\right\}\,.
\endaligned
\tag7.15
$$
By making the simplest choice $f(z) = 1 - |z|^2$ for the equation of the boundary, we
see that
$\cc(Rj_*\CC_{D_-})$ is a cylinder with base $D_-$. There is a more ``sophisticated"
choice for the equation of the boundary:
$$
  f(z) = \left(\frac{1+ |z|^2}{1-|z|^2}\right)^{4\lambda}\,.
\tag7.16
$$
For this boundary equation we get
$$
\mu_\lambda\left\{ \frac{df_x}{f(x)} \mid x \in D_-\right\}\  = \
\Omega_\lambda(D_-)\,,
\tag7.17
$$
where, we recall, $\Omega_\lambda(D_-)$ denotes the coadjoint $\GR$-orbit in
$i\fg^*_\Bbb R$ determined by $\lambda$ and $D_-$. Letting the parameter $s$
vary between 0 and 1 establishes a homology between $\Omega_\lambda(D_-)$ and
$\mu_\lambda(\cc(Rj_*\CC_{D_+}))$. Thus, our character formula (7.13) agrees with
Rossmann's formula (7.3) for discrete series representations of $SL(2,\Bbb R)$.

The above argument  can be generalized to  discrete series representations for any
group $\GR$.  To this end, let us assume that $\GR$ has a compact Cartan. Pick an  open
$\GR$-orbit 
$S\subset X$ and choose $\lambda \in \rho +\Lambda$ antidominant. Recall that the
discrete series representation attached to $(S,\lambda)$ is given by
$(-1)^s M(Rj_*\CC_{S})$, where $j:  S \hookrightarrow X$ denotes the inclusion.
In analogy with (7.16), we make the following choice for the defining
equation of
$\partial S$:
$$
f\ =\ \frac{\text{$\GR$-invariant metric on $\bold
L_\lambda$}}{\text{$U_\RR$-invariant  metric on $\bold L_\lambda$}}\,;
\tag7.18
$$
here $\bold L_\lambda$ denotes the line bundle on $X$ associated to $\lambda$. Then,
just as in the case $\GR=SL(2,\Bbb R)$, 
$$
\mu_\lambda\left\{ \frac{df_x}{f(x)} \mid x \in S\right\}\  = \
\Omega_\lambda(S)\,.
\tag7.19
$$
Therefore, for discrete series, the character formula (7.13) agrees with Rossmann's
formula.  This is a crucial step in proving the  formulas (7.13) and (7.14). For
details, see
\cite{SV2}.

\Lecture8{Microlocalization of Matsuki = Sekiguchi}			

In this final lecture we explain the relationship between the Matsuki and the
Sekiguchi correspondences. This relationship is provided by geometry. It is a
crucial step in the proof of the Barbasch-Vogan conjecture \cite{SV3} which we will
briefly explain at the end. We will continue to use the notation of the previous
lectures.  Recall that  the Matsuki correspondence provides a bijection between the
$\GR$- and $K$-orbits on $X$:
$$
\GR\backslash X \ \longleftrightarrow \ K\backslash X\ ,
\tag8.1
$$
where a $\GR$-orbit $\O'$ corresponds to $K$-orbit $\O$ if and only if $\O\cap\O'$
is non-empty and compact (in which case  $\O\cap\O'$ constitutes a $K_\RR$-orbit). Let
$\Cal N$ be the nilpotent cone in $\fg$. The Kostant-Sekiguchi correspondence
provides  a bijection between  nilpotent $\GR$- and
$K$-orbits: 
$$
  \GR \backslash  \Cal N \cap i\gg_\RR\ \longleftrightarrow\ K\backslash  \Cal N
\cap \fp
\tag8.2
$$
Here a $\GR$-orbit $\O'$ corresponds to a $K$-orbit $\O$ if and only if
there exists a Lie algebra homomorphism $\frak s\frak l_2(\CC) \rightarrow \gg$
which is defined over $\RR$ and which commutes with the Cartan involution,
such that 
$$
  j \left(\matrix 0 & i \\  0 & 0 \endmatrix\right) \in \O'\ , \qquad
  j \left(\matrix 1 & i \\  i & -1 \endmatrix\right) \in \O\ .
\tag8.3
$$
Let us recall, from lecture 5, the Matsuki correspondence for sheaves:
$$
\gamma\,:\,\od_K(X)_{-\lambda}\  @>{\ \ \sim \ \ }>> \  \od_\GR(X)_{-\lambda}\,.
\tag8.4
$$
It satisfies the following property: 
$$
  \gamma(Rj_*\CC_{\O}) = Rj'_!\CC_{\O'} [2 \operatorname{codim}_\CC \O]\,;
\tag8.5
$$
here $j:\Cal O \hookrightarrow X$ and $j':\Cal O' \hookrightarrow X$ denote the
inclusions of a $K$-orbit $\Cal O$ and a $\GR$-orbit $\Cal O'$ which are related under
the Matsuki correspondence. Also, recall that
$\gamma=\Gamma^{-1}$ is defined exactly in the same way as $\Gamma$ except that the
roles of $K$ and $\GR$ are switched and all the *'s are replaced by !'s. Next, we
consider the following commutative diagram:
$$
  \CD
  \od_K(X)_{-\lambda} @>\gamma>>  \od_\GR(X)_{-\lambda} \\
   @VV{\cc}V @VV{\cc}V \\
   \oh^{inf}_{2n}(T^*_K X,\ZZ)  @>{\Phi=\cc(\gamma)}>> \oh^{inf}_{2n}(T^*_\GR
X,\ZZ)\,,
   \endCD
\tag8.6
$$
where the map $\Phi=\occ(\gamma)$ is the effect of the functor $\gamma$ on
characteristic cycles. From the properties (a-d) of the characteristic cycles one can
conclude that apriori that $\Phi=\cc(\gamma)$ exists. We will now give an explicit
formula for
$\Phi$. To this end, let us define, for $s>0$, an automorphism $F_s: T^*X \rightarrow
T^*X$ by the following formula:
$$
  F_s(x,\xi)\ = \ (\, \exp(s^{-1}(\Re \xi))\, x\, ,\,  
  \Ad(\exp(s^{-1}(\Re \xi))\,\xi\, )\ , 
\tag8.7
$$
where $(x,\xi)\in  T^*X, \xi\in T_x^*X \simeq \fn_x$, and
$\Re: \gg \rightarrow \gg_\RR$ associates to an element in $\fg$ its real part in
$\gr$. The automorphisms $F_s$ preserve the real(!)  symplectic form on
$T^*X$. The map
$\Phi$ is given as a limit of these symplectomorphisms:
$$
  \Phi(C) \ = \ \lim_{s \rightarrow 0^{+}} (F_s)_*(C)\ .
\tag8.8
$$
\remark{Remark}
It is not obvious that the limit exists. To give meaning to it, one constructs the
chain $\tilde C = \{(s,\zeta)\in \Bbb R\times \fg \mid s>0\,,\  \zeta\in (F_s)_*(C)\}$
and sets $ \lim_{s \rightarrow 0^{+}} (F_s)_*(C) = - \partial\tilde C$. For this to
make sense, the support $|\tilde C|$ should be triangulable. This is by no means
obvious. It follows from the fact that the set $|\tilde C|$ belongs to the
analytic geometric category $\Cal C$ coming from the o-minimal structure $\RR_{an,
\exp}$. For a beautiful exposition of analytic geometric structures see \cite{DM}. 
The analytic geometric categories satisfy (essentially) the same good properties as
semi-algebraic (and subanalytic) sets, and we could have worked in this more
general context in lectures 2 and 6. The reason that one needs to pass to  the
analytic geometric category $\Cal C$ in Lie theory is that the exponential function
$\exp:\Bbb R \to \Bbb R$ is not real analytic at infinity. However, it is an
allowable function in the category $\Cal C$.
\endremark

We will briefly explain the idea behind the verification of the  formula (8.8), for
the details, see \cite{SV3}. To begin with, let us recall the definition of
$\gamma$. Using the notation in the diagram:
$$
X @<a<<\GR\times X @>q>> \GR/\KR \times X @>p>> X
\tag8.9
$$
the functor $\gamma$ is given by $\gamma(\F) = Rp_! \tilde{\F}$ where
$\tilde{\F}$ is the object such that $q^{!}\tilde{\F}=a^!\F$. As the maps $a$ and
$q$ are submersions, it  is easy to express $\occ(\tilde\cf)$ in terms of
$\occ(\cf)$.  The difficulty lies in computing $\cc(Rp_!\tilde{\F})$ in terms of
$\occ(\tilde\cf)$, in particular, because the map $p$ is not proper. We will
compactify the map $p$ by embedding  the symmetric space $\GR/\KR$ inside a compact
real algebraic manifold and taking its closure. At this point it does not matter
which compactification of the symmetric space we choose. This gives us the following
diagram:
$$
\CD 
\GR\times X @>q>> \GR/\KR \times X @>j>> \overline{\GR/\KR}\times X \\
@VaVV             @VVpV               @VV\bar{p}V \\
X @. X @= X\,.
\endCD
\tag8.10
$$ 
From this diagram we conclude that
$$
  Rp_!\tilde{\F} = R\bar{p}_* j_! \tilde{\F}\,,
\tag8.11
$$
using the fact that $\bar{p}_!=\bar{p}_*$ since $\bar{p}$ is proper. To get  a 
formula for $\occ(Rp_!\tilde\cf)=\occ(R\bar p_*j_!\tilde \cf)$, we begin by applying
(6.11) to $R\bar p_*$. The top row of (6.13) becomes, in this case,
$$
  T^*(\overline{\GR/\KR} \times X) @<{\ \ d\bar p^*\ }<< 
  \GR/\KR\times T^*X @>{\ \tau\ }>> T^*X\,,
\tag8.12
$$
and $d\bar p^*$ is the canonical embedding. Then, by (6.11), 
$$
\occ(R\bar p_*j_!\tilde \cf)=\bar{p}_*(\occ(j_!\tilde \cf)) =
\tau_*([\overline{\GR/\KR} \times X]
. \occ(j_!\tilde \cf))\,,
\tag8.13
$$
where, $[\overline{\GR/\KR} \times X]. \occ(j_!\tilde \cf)$ is the intersection
product of cycles.  To compute $\occ(j_!\tilde \cf)$, we choose a defining equation
$f$ for the boundary of $\GR/\KR$ and  use  a variant of
the defining property (d) of the characteristic cycle map. Denoting 
$C=\cc(\tilde \F)$, we get:
$$
\occ(j_!\tilde\cf)\ = \ j_! C\  = \ \lim_{s\rightarrow 0^+} (C -s \,d\log{f})\,.
\tag8.14
$$ 
Combining (8.13) and (8.14) gives:
$$
  \cc(Rp_!\tilde{\F}) 
  \ = \ \tau_*([\overline{\GR/\KR} \times X] . (\lim_{s\rightarrow 0^+} (C -s
\,d\log{f})))\,.
\tag8.15
$$
We now rewrite this formula as
$$
  \cc(Rp_!\tilde{\F}) 
  \ = \ \lim_{s\rightarrow 0^+}\tau_*([\overline{\GR/\KR} \times X] .  (C -s
\,d\log{f}))\,.
\tag8.16
$$
If $f$ is chosen appropriately then the intersection $(\overline{\GR/\KR} \times
X)\cap (C -s\,d\log{f})$ is transverse and 
the intersection product in (8.16)  can be replaced by an ordinary
intersection\footnote{We are ignoring all issues of orientation.}. This
observation is the crux of the computation.

To evaluate the right hand side of (8.16), we choose as $\overline{\GR/\KR}$ the
one-point compactification of $\GR/\KR \simeq\frak{p}_\RR$:
$$
  \overline{\GR/\KR}\ \simeq\ \overline{\frak{p}_\RR}\ =\  \frak{p}_\RR \cup
\{\infty\}
\tag8.17
$$
Furthermore, we choose  the defining equation $f: \overline{\frak{p}_\RR}
\rightarrow \RR$ as follows:
$$
  f(\zeta) =  \cases e^{-\frac{1}{2}B(\zeta,\zeta)} &\text{if $\zeta
\in\frak{p}_\RR$} \\
0 &\text{if $\zeta=\infty$}\,.\endcases
\tag8.18
$$
Here $B$ denotes the Killing form. Then  $d \log f = - \id $ on
$\frak{p}_\RR$ and a relatively easy computation gives the formula for $\Phi$.

\remark{Remark} The function $f$ is not real analytic at infinity. This forces us
outside the semi-algebraic and subanalytic categories and into the analytic
geometric category $\Cal C$ coming from the o-minimal structure $\Bbb
R_{\text{an,exp}}$, as was explained in the previous remark. 
\endremark

We will now relate the map $\Phi$ to the Kostant-Sekiguchi correspondence. To do so, we
identify $\fg$ with $\fg^*$ via the Killing form and let $\Cal N\subset \fg$ denote
the nilpotent cone. Then the moment map $\mu:\ct \to \Cal N$ and furthermore
$$
\mu^{-1}(i\gg_\RR)\  = \ T^*_\GR X\qquad\text{and} \qquad
\mu^{-1}(\frak p) = T^*_K X\,.
\tag8.19
$$ 
Let us fix a $G$-orbit $\Cal O \subset \Cal N$, write $i\fg_\RR \cap \Cal O$
as a union of $\GR$-orbits:
$$
  i\fg_\RR \cap \Cal O = \O'_1 \cup \cdots \cup \O'_k\,, 
\tag8.20a
$$
and, similarly, $\frak p \cap \Cal O$ as a union of $K$-orbits:
$$
  \frak p \cap \Cal O = \O_1 \cup \cdots \cup \O_k\,.
\tag8.20b
$$ 
We enumerate the orbits so that $\O_i$ and $\O'_i$ correspond to each other under 
Sekiguchi. Let us consider the complex Lagrangian cycle $[\mu^{-1}(\Cal
O_i)]$ on
$T^*X$. It is, by definition, supported on $T^*_KX$. Clearly the symplectomorphisms
$F_s$ of (8.7) map $ \mu^{-1} (\O) \rightarrow \mu^{-1} (\O)$ and hence $\Phi$, as the
limit of the $F_s$, maps $\mu^{-1} (\O)$ to its closure $\overline{\mu^{-1}(\O)}$\,.
Thus, we can write
$$
  \Phi([\mu^{-1}(\O_i)]) = \sum n_j [\mu^{-1} \O'_j] + \text{lower order terms}\,,
\tag8.21
$$
where by lower order terms we mean chains which lie over $\partial \bar\Cal O =
\bar\Cal O-\Cal O$. Here the $n_j\in\Bbb Z$ and we note that the $[\mu^{-1} \O'_j]$
are chains, not necessarily cycles, as they can have boundary in $\mu^{-1}(\partial
\bar\Cal O)$.

\proclaim{Theorem} The map $\Phi$ induces the Sekiguchi correspondence on the
nilpotent orbits, i.e., 
$$
  \Phi([\mu^{-1}\O_i]) = [\mu^{-1}\O_i'] + l.o.t\,.
\tag8.22
$$
\endproclaim
We do not know of a simple argument for (8.22). The proof is contained in
\cite{SV3,SV4}. As a corollary of this result we see that the Sekiguchi
 correspondence is given by
$$
  c\mapsto 
\lim_{s\rightarrow 0^+} \{\Ad(\exp(s^{-1} \Re \zeta) \zeta  \mid \zeta \in c\} =
\lim_{s\rightarrow 0^+} \{s\Ad(\exp(\Re \zeta) \zeta  \mid \zeta \in c\}\,,
\tag8.23
$$
where $c$ stands for one of the $K$-orbits in (8.20b).

Finally, we explain, very briefly, how the above theorem enters the proof of the
Barbasch-Vogan conjecture. For details, see \cite{SV3}. Putting together (5.4),
(8.6), and (8.22) we get the commutative diagram:
$$
  \CD
  \{ \text{HC-modules}\}_{\chi_\lambda} @>>>
\{\text{$\GR$-representations}\}_{\chi_\lambda}
\\
  @VVV @VVV \\
  \od_K(X)_{-\lambda} @>\gamma>>  \od_\GR(X)_{-\lambda} \\
   @V{\cc}VV @VV{\cc}V \\
  \oh^{inf}_{2n}(T^*_K X,\ZZ)  @>{\Phi}>> \oh^{inf}_{2n}(T^*_\GR X,\ZZ) \\
   @V{\operatorname{gr}(\mu)_\lambda}VV @VV{\operatorname{gr}(\mu)_\lambda}V \\
  \text{nilpotent orbits in $\Cal N \cap \frak p$}
   @>{\ \ \text{Sekiguchi}\ \ }>>
   \text{nilpotent orbits in $\Cal N \cap \frak i\fg_\RR$}
\endCD
\tag8.24
$$
A few remarks are in order. First, we have turned around the arrows in (5.4). As
the functors $M$ and $\alpha$ are not invertible, one interprets the top square as
a consistent choice of representatives in $\od_K(X)_{-\lambda}$ and $
\od_\GR(X)_{-\lambda}$ for representations. Furthermore, both the vertical arrows
${\operatorname{gr}(\mu)_\lambda}$ stand for the operation of taking the leading term of
the result of integration of a cycle along the fiber of $\mu$ against the form
$e^\lambda$.

Let us consider an irreducible representation $(\pi,V)$ whose
associated Harish-Chandra module we denote by $M$. By a result of Chang
\cite{C}, the left vertical column amounts to the associated cycle construction:  
$$
\operatorname{Ass}(M) \ = \ \sum a_j [\Cal O_j]\,, \qquad \text{with} \ a_j\in\Bbb
Z_{\geq 0}\,.
\tag8.25
$$ 
For the associated cycle  construction, see the lectures of Vogan.
From the character formula (7.14) we can
conclude that the right hand column amounts to associating to $V$ the Fourier
transform of the leading term of the asymptotic expansion of the Lie algebra
character of the representation $(\pi,V)$. This gives the wave front cycle
$$
\operatorname{WF}(V) \ = \ \sum b_j [\Cal O'_j]\,, \qquad \text{with} \ b_j\in\Bbb
C\,.\tag8.26
$$
The invariant $\operatorname{Ass}(M)$ is purely algebraic whereas the invariant
$\operatorname{WF}(V)$, introduced in \cite{BV}, is analytic. The commutative
diagram (8.24) implies that these two invariants coincide under the Sekiguchi
correspondence, i.e.,
$$
a_j\ =\ b_j \qquad \text{if $\ \ \O_j \leftrightarrow \O'_j\ \ $ under Sekiguchi}\,.
\tag8.27
$$
This statement is usually referred to as the Barbasch-Vogan conjecture. Note that
it implies, in particular, that the $b_j$ are non-negative integers.

\def\Appendix#1{
   \topmatter
   \lecturelabel{Appendix}
   \lecture{{}\!\!}
   \lecturename #1\endlecturename
   \endtopmatter
}
\Appendix{Homological algebra\\by Markus Hunziker}

\noindent
Let $\frak A$ be an abelian category. 
Typical examples of abelian categories are
the category 
of (left) modules over an arbitrary ring $R$ with unit,
the category of sheaves of $\CC$-vector spaces on a topological space $X$, and
the category of $\CC$-constructible 
sheaves on a semi-algebraic set $X$.

\head The category of complexes $\operatorname{C}(\frak A)$\endhead

\definition{Definition}
Recall that a ({\it cochain\/}) {\it complex\/} of objects in $\frak A$ is a 
sequence of objects $A^i$, $i\in\ZZ$, together with morphisms $d^{@,@,i}: A^i
\rightarrow A^{i+1}$,
$$
A^\cdot \ = \ (\, \cdots @>{\ \ }>> A^i @>{\,d^{i}\,}>> A^{i+1}
@>{d^{@,@,i+1}\!}>>
  A^{i+2} @>{\ \ }>> \cdots\,) \, , \quad
 \text{such that $\ d^{@,@,i+1}\circ d^{@,@,i} =0$}\ 
$$
for all $i$.
The morphisms $d^{@,@,i}: A^i \rightarrow A^{i+1}$ are called the 
{\it differentials\/} of the complex.
A {\it morphism of complexes\/}  $f: A^\cdot \rightarrow B^{@,@,\cdot}$ 
is a sequence of morphisms $f^i: A^i  \rightarrow B^{@,@,i}$
which commute with the differentials in the sense that
$ d^{@,@,i}_{\!B}\circ f^i =  f^{i+1}\circ d^{@,@,i}_{\!A} $ for all $i$.
Thus we obtain a category $\operatorname C(\frak A)$,
which is  abelian.
We identify $\frak A$ with the full subcategory of $\operatorname C(\frak A)$
consisting of complexes $A^\cdot$ such that $A^i = 0\ $ for $i \not =0\,$.
Later we will also need the full subcategory $\operatorname{C}^+\!(\frak A)$ 
of bounded below complexes. 
\enddefinition

\definition{Shift functors}
For every integer $k$, we define a functor 
$[\,k\,]: \operatorname C(\frak A) @>>> \operatorname C(\frak A)$ as follows.
If $A^\cdot $ is a complex then $A^\cdot[@,@,k@,@,]$ is the complex given by
$$
  A^i[@,@,k @,@,]\ =\ A^{k+i}\ , \quad \quad  
  d^i_{\!A[@,@,k @,@,]}\ =\  (-1)^k\, d_{\! A}^{k+i} \ .
$$
If $f: A^\cdot \rightarrow B^{@,@,\cdot}$  is a morphism of complexes then
$f[@,@,k@,@,]:  A^\cdot[@,@,k @,@,] @>>> B^{@,@,\cdot}[@,@,k @,@,]$ is given
by $f^i[@,@, k@,@,]  = f^{k+i}$.
The functor $[@,@,k @,@,]$ is called the {\it shift functor\/} of degree $k$.
\enddefinition

\definition{Cohomology and quasi-isomorphisms}
The {\it $i$-th cohomology object\/} of a complex $A^\cdot$ is the object
$H^i(A^\cdot)  = \ker d^{@,@,i}/\im d^{@,@,i-1}$ which is well-defined since
$\frak A$ is an abelian category.
A morphism of complexes $f: A^\cdot \rightarrow B^{@,@,\cdot}$ 
induces morphisms $ H^i(f) : H^i(A^\cdot) \rightarrow H^i(B^{@,@,\cdot})$
between cohomology objects. 
If all the $ H^i(f) : H^i(A^\cdot) \rightarrow H^i(B^{@,@,\cdot})$
are isomorphisms then we say that $f$ is a {\it quasi-isomorphism\/} and
we write $f: A^\cdot \qis B^{@,@,\cdot}$.
\enddefinition

\definition{Example} Let $A$ be an object of $\frak A$ and let
$\ 0 @>>> A @>>> E^{0} @>>> E^{1} @>>> \cdots\ $ be a resolution
of $A$. Then we have a quasi-isomorphism $A \qis E^{@,@,\cdot}$.
\enddefinition

\newpage 
\head The homotopy category $\operatorname{K}(\frak A)$ \endhead

\definition{Homotopy}
Two morphisms $f,g: A^\cdot \rightarrow B^{@,@,\cdot}$ in
$\operatorname{C}(\frak A)$ are  called {\it homotopic\/} if there 
is a sequence of morphisms
$k^i:A^i \rightarrow B^{i-1}$ in $\frak A$
such that
$$
  f^i - g^i\ =  \ d^{\,i-1}_{\!B} \circ k^i \,+\, k^{i+1} \circ d_{\!A}^i\ .
$$
If $f$ and $g$ are homotopic then they induce the same morphism
$H^i(A^\cdot) \rightarrow H^i(B^{@,@,\cdot})$ on the cohomology
objects for all $i$.
\enddefinition

\definition{Definition}
The category $\operatorname{K}(\frak A)$ has as objects complexes 
of objects in $\frak A$ and as morphisms homotopy equivalence classes 
of morphisms in $\operatorname{C}(\frak A)$. Similarly, we obtain
a category $\operatorname{K}^+\!(\frak A)$ from the category
$\operatorname{C}^+\!(\frak A)$  of bounded below complexes.
\enddefinition

\head Triangles and long exact sequences \endhead
The notion of a short exact sequence is  not well-defined
in $\ok(\frak A)$, which is {\it not\/} an abelian category.
The substitutes for short exact sequences are so-called distinguished triangles.
They generate canonically long exact cohomology sequences.

\definition{Definition}  A {\it triangle\/} in $\ok(\frak A)$ is a diagram
$ A^\cdot  @> u >>  B^{@,@,\cdot}  @> v >>  C^{@,@,@,\cdot}
  \ @> w >> A^\cdot[\,1\,]$\ .
Often a triangle is written in the mnemonic form
$$
  \matrix 
 & A^\cdot \ @>{\quad\! u\quad\!}>> \ B^{@,@,\cdot} & \\ \vspace{1 pt}
 & \!{}_{w}\!\!\nwarrow^{\!\!\!+1} \,\,\quad \!\swarrow \!\!\!{}_{v} \  & \
 \\ \vspace{1 pt}
 &        \!  C^{@,@,@,\cdot}               &
 \endmatrix  \!\!\!,
$$
whence the name.
A {\it morphism of triangles\/} is given by a commutative diagram:
$$
  \CD
   A^\cdot  @> u >> B^{@,@,\cdot} @> v>> C^{@,@,@,\cdot}
   @>\ w\ >> A^\cdot[\,1\,] \\
  @VV{f}V   @VV{g}V   @VV{h}V @VV{f[\,1\,]}V \\
   \,{A'}^\cdot\! @> u' >> \,{B'}^{\cdot}\!  @> v' >> 
   \,{C'}^{\cdot}\!
  @> w' >> {A'}^\cdot[\,1\,]
  \endCD
$$
\enddefinition

\definition{Mapping cones and distinguished triangles} 
Let $u : A^\cdot \longrightarrow B^{@,@,\cdot}$ be a morphism in 
$\operatorname{C}(\frak A)$. 
The {\it mapping cone} of $u$ is the complex $C(u)^\cdot$
which is defined by
$$
   \,C(u)^i\ = \  A^{i+1}\,\oplus\,  B^{@,@,i}
\quad, \quad \quad
  \quad d_{C(u)}^i \ = 
  \left[{\matrix
  \!-d_A^{@,@,i+1} & 0 \\ \vspace{4 pt}
  \,\,\,u^{i+1}  & d_B^{@,@,i}\,
  \endmatrix}\right]\quad .
$$
There is a canonical exact sequence 
$\ 0 @>>> B^{@,@,\cdot} @>v>> C(u)^\cdot @>w>> A^\cdot[\,1\,] @>>> 0\ $ in
$\operatorname{C}(\frak A)$, given by $v:b\mapsto (0,b)$ and $w: (a,b) \mapsto a$.
The triangle 
$$
   A^\cdot \ @>\ u\ >>\  B^{@,@,\cdot} \ @>\ v\ >>\  C(u)^\cdot
  \ @>\ w\ >>\ A^\cdot[\,1\,]\ 
$$
is called the {\it standard triangle\/} associated to the mapping cone $C(u)^\cdot$.
A {\it distinguished triangle} in $\operatorname{K}(X)$ is a triangle
which is isomorphic to a standard one. 
\enddefinition

\remark{Remark}
If two morphisms $u, u':  A^\cdot \rightarrow B^{@,@,\cdot}$ in
$\operatorname{C}(\frak A)$ are homotopic then the mapping cones 
$C(u)^\cdot$ and $C(u')^\cdot$ are  isomorphic in $\ok(\frak A)$
and also their associated standard triangles are isomorphic.
This isomorphism is not unique in general.
\endremark

\definition{Long exact sequences} 
If 
$\ A^\cdot @>u>>  B^{@,@,\cdot}@>v>>C^{@,@,@,\cdot}
   @>w>>A^\cdot[\,1\,]\ $
is a triangle in $\ok(\frak A)$, 
then the morphisms $u$, $v$, and $w$ induce canonically a sequence
$$
  \cdots \ @>>>\ H^i(A^\cdot)\ 
  @>>>\ H^{i}(B^{@,@,\cdot})\ @>>>\
   H^i(C^{@,@,@,\cdot})\
  @>>>\
  H^{i+1}(A^\cdot)\ @>{\ \ }>>\ \cdots\, \ .
$$
If the triangle is distinguished this sequence is exact. 
\enddefinition

\remark{Remark}
Let
$\, 0  @>>> A^\cdot @>u >>  B^{@,@,\cdot} 
  @>v>> C^{@,@,@,\cdot} @>>> 0\, 
$
be a short exact sequence in $\operatorname{C}(\frak A)$.
There is a canonical map $h: C(u)^\cdot \rightarrow C^{@,@,@,\cdot}$ given
by $h^i: A^{i+1}\oplus B^{@,@,i}@>>> C^i,\ (a,b)\mapsto v^i(b)$.
One can show that $h$ is a quasi-isomorphism.
\endremark

%

%
%

\head The derived category $\od(\frak A)$\endhead
The derived category $\od(\frak A)$ is obtained from $\ok(\frak A)$
by ``localization'' at the multiplicative set of quasi-isomorphisms.
It comes together with a natural functor $Q: \ok(\frak A) \rightarrow \od(\frak A)$ 
which sends quasi-isomorphisms to isomorphisms. Similarly, a derived 
category $\od^+\!(\frak A)$ is obtained from $\ok^+\!(\frak A)$.

\definition{Definition} 
The objects of $\od(\frak A)$ are again just complexes of objects in $\frak A$.
If $A^\cdot$ and $ B^{@,@,\cdot}$ are two objects in  $\od(\frak A)$ then
a morphism from $A^\cdot$ to $B^\cdot$ is defined as an equivalence classes of
diagrams in $\ok(\frak A)$ of the form $A^\cdot\leftqis C^{@,@,@,\cdot} \rightarrow
B^{@,@,\cdot}$. Here the diagram $A^\cdot\leftqis C^{@,@,@,\cdot}
\rightarrow B^{@,@,\cdot}$ is equivalent to the diagram  
$A^\cdot\leftqis \tilde C^{@,@,@,\cdot} \rightarrow B^{@,@,\cdot}$ if there
exists a commutative diagram in $\ok(\frak A)$ of the form
$$
  \matrix 
 &          \, C^{@,@,@,\cdot}               & \\ \vspace{1 pt}
 & \swarrow \ \ \uparrow \ \ \searrow & \\ \vspace{5 pt}
 & A^\cdot\!  @<{\,\,{}_{qis}}<< \, D^{@,@,\cdot}  @>{\ \ }>> B^{@,@,\cdot} &\\ 
 \vspace{4 pt}
 & \nwarrow \ \ \downarrow \ \ \nearrow & \\ \vspace{3 pt}
 &         \, \tilde C^{@,@,@,\cdot}   &
 \endmatrix
$$
One can view a diagram $A^\cdot @<s<< C^{@,@,@,\cdot} @>u>> B^{@,@,\cdot}$ 
as a fraction $u/s$.
The composition of morphisms in $\od(\frak A)$ is defined as follows.
Let $A^\cdot \leftqis D^{@,@,\cdot} @>>> B^{@,@,\cdot}$ and 
$B^{@,@,\cdot} \leftqis E^{@,@,\cdot} @>>> C^{@,@,@,\cdot}$
be two diagrams in $\operatorname{K}(\frak A)$
representing two morphisms in $\od(\frak A)$.
One can show that there always exists a diagram $D^{@,@,\cdot} \leftqis
F^{@,@,\cdot} @>>> E^{@,@,\cdot}$ such that the following diagram commutes in
$\operatorname{K}(\frak A)$:
$$
  \matrix
  &     \,F^{@,@,\cdot} & \\
  &  \swarrow \ \quad \ \  \searrow & \\ \vspace{3 pt}
  & \, D^{@,@,\cdot} \ \quad \quad \quad \quad \ E^{@,@,\cdot}  & \\
  &   \swarrow \ \quad \ \  \searrow \ \quad \ \swarrow \ \quad \ \ 
\searrow &
\\
\vspace{3 pt}
  &  A^\cdot \ \quad \quad \quad \quad \ \, B^{@,@,\cdot} \ \quad \quad \quad \quad
\ C^{@,@,@,\cdot}   
\endmatrix
$$
\vskip 5 pt
\noindent
The diagram $A^\cdot \leftqis F^{@,@,\cdot} @>>> C^{@,@,@,\cdot}$ then defines
the composition of the given morphisms in $\od(\frak A)$.
Every morphism $A^\cdot @>>>
B^{@,@,\cdot}$ in $\operatorname{K}(\frak A)$ induces a morphism in $\od(\frak A)$ via
the diagram $A @<{}_=<< A @>>>B$. This defines the functor $Q: \ok(\frak A) @>>>
\od(\frak A)\,$.
\enddefinition

\remark{Remark}
A morphism $u: A^\cdot \rightarrow B^{@,@,\cdot}$ in $\ok(\frak A)$
becomes the zero-map in $\od(\frak A)$ iff there is a quasi-isomorphism 
$s: B^{@,@,\cdot} \rightarrow A^\cdot$ such that 
$s \circ u = 0$ in  $\ok(\frak A)$.
\endremark

\definition{Hyperext and homological dimension}
Let $A^\cdot$ and $B^\cdot$ be two complexes considered as objects in the 
derived category $\od(\frak A)$. Then we define
the {\it $k$-th hyperext\/} as the abelian group
$$
  \Ext^k(A^\cdot,B^\cdot)\ = \ \Hom_{\od(\frak A)}(A^\cdot,B^\cdot[k])\ .
$$
If $A$ and $B$ are objects in $\frak A$, which we may consider as objects
in $\od(\frak A)$ concentrated in degree $0$, then $\Ext^k(A,B)$ coincides
with the usual $\Ext$. This is a result due to Yoneda.

We say that $\frak A$ has {\it homological dimension\/} $\leq n\,$ if
$\Ext^k(A,B)= 0$ for all $k>n$ and for any objects $A$
and $B$ in $\frak A$. If $\frak A$ has homological dimension $\leq 1$
then one can show that for any complex $A^\cdot$ which is bounded from above 
and below we have an isomorphism in the derived category
$$
  A^\cdot\ \simeq \ \bigoplus_k H^k(A^\cdot)[-k@,@,]\ .
$$
This holds for example if $\frak A$ is the abelian category of vector spaces 
over a field.
\enddefinition

%
%

\head Derived functors \endhead
Let $F: \frak A @>>> \frak B$ be an additive morphism between abelian 
categories. The functor extends to a functor 
$\ok^+\!(\frak A) @>>> \ok^+\!(\frak B)$. However, this functor does not
send quasi-isomorphisms to quasi-isomorphisms in general.
If we assume that the functor $F$ is left exact then under 
suitable hypotheses (for example if $\frak A$ has enough injectives)
there exists a derived functor
$RF: \od^+\!(\frak A) \rightarrow \od^+\!(\frak B)$ which is close to $F$
in the sense that if $A$ is an object in $\frak A$ we have a natural isomorphism
$F(A)= H^0(RF(A))$.

\definition{Injective resolutions}
An {\it injective resolution\/} 
is a quasi-isomorphism
$A^\cdot \qis I^{@,@,\cdot}$ such that $I^i$ is an injective object of
$\frak A$ for all $i$. If $\frak A$ has enough injectives then injective
resolutions always exist. 
In $\ok(X)$ injective resolutions are also unique as follows.
Suppose $f: A^\cdot @>>> B^{@,@,\cdot} $ is a morphism in $\ok(\frak A)$.
Let $u: A^\cdot \qis E^{@,@,\cdot}$ be any quasi-isomorphism and let 
$v: B^{@,@,\cdot} \qis I^{@,@,\cdot}$ be an injective resolution.
Then there exists a unique morphism $g: E^{@,@,\cdot} @>>> I^{@,@,\cdot}$
such that $v\circ f = g\circ u$.
\enddefinition

\proclaim{Theorem}
Assume $\frak A$ has enough injectives.
Let $\frak J$ be the full category of $\frak A$ of injective objects.
Then the natural functor $Q:\ok^+\!(\frak A) @>>> \od^+\!(\frak A)$
induces an equivalence of categories
$$
   \ok^+\!(\frak J) \ \simeq \ \od^+(\frak A)\ .
$$
\endproclaim

For a proof of this theorem see, for example, \cite{KSa, Prop. 1.7.10}


\definition{Definition}
Let $F: \frak A \rightarrow \frak B$ be a left exact, additive functor between abelian
categories and assume that $\frak A$ has enough injectives.
Then the {\it right derived functor\/} of $F$ is the functor 
$RF: \od^+\!(\frak A) \rightarrow \od^+\!(\frak B)$ given by
$$
 RF(A^\cdot)\ =\  F( I^{@,@,\cdot})\ ,
$$
where $I^{@,@,\cdot}$ is any injective resolution of $A^\cdot$.
Note that $RF$ is well-defined by the remarks above.
The $i$-{\it th derived functor\/} of $F$ is the functor 
$R^i\!F: \od^+\!(\frak A) \rightarrow \frak B$
given by
$ R^i\!F(A^\cdot) = H^i(F(I^{@,@,@,\cdot}))$. 
\enddefinition

\remark{Remark} More generally, one can define define the derived 
functor $RF$ for a functor $F: \ok^+\!(\frak A) @>>> \ok^+\!(\frak B)$
of triangulated categories, i.e., a functor which commutes with
the shift functor $[\,1\,]$, and transforms distinguished triangles 
into distinguished triangles.
\endremark

\remark{Example}
Fix an object $A^\cdot$ of $\ok^+\!(\frak A)$.
Then the functor $\Hom(A^\cdot,-): \ok^+(\frak A) \rightarrow \ok^+(\frak A)$ 
is a functor of triangulated  categories and we may compute the 
derived functor $R\!\Hom(A^\cdot,-)$ using injective resolutions as above.
\endremark

\proclaim{Theorem (Yoneda)}
Assume that $\frak A$ has enough injectives.
Then 
$$
  \Ext^k(A^\cdot ,B^\cdot) \ =\ R^k\!\Hom(A^\cdot ,B^\cdot)\ .
$$ 
In particular, for two objects $A$ and $B$ of $\frak A$,
$\Ext^k(A ,B ) $ is the usual $\Ext$.
\endproclaim
 

\remark{$F$-injective resolutions}
To compute the derived functor of $F$ it not necessary to consider
injective resolutions.
A full additive subcategory $\frak J$ of $\frak A$  is called 
$F$-{\it injective} if the following conditions are satisfied:
\vskip 5 pt
\roster
\item"(i)" Every object of $\frak A$ is isomorphic to a subobject of an
object of $\frak J$.
\item"(ii)" If $\ 0 @>>> A' @>>> A @>>> A'' @>>> 0 \ $ is an exact sequence
in $\frak A$, and if $A'$ and $A$ are objects of $\frak J\,$, then $A''$ 
is also an object of $\frak J\,$.
\item"(iii)" If $\ 0 @>>> A' @>>> A @>>> A'' @>>> 0\ $ is an exact sequence in $\frak
A$, and if $A'$, $A$, $A''$ are objects in $\frak J$, then the sequence $\,0 @>>>
F(A') @>>> F(A) @>>> F(A'') @>>> 0\,$ is exact.
\endroster
\vskip 5 pt
\noindent
If $\frak J$ is any $F$-injective subcategory of $\frak A$ then we may
compute
$RF:
\od^+\!(\frak A) \rightarrow \od^+\!(\frak B)\,$ as above by replacing injective
resolutions with $F$-injective resolutions. 
This is very useful to compute derived functors in practice by choosing a 
convenient $F$-injective category. 
\endremark

\proclaim{Theorem}
Let $F:\frak A @>>> \frak A'$ and $F': \frak A' @>>> \frak A''$ be two 
left exact additive functors between abelian categories. Assume that
there exists 
an $F$-injective subcategory $\frak J$ of $\frak A$, and 
a  $F'$-injective subcategory $\frak J'$ of $\frak A'$ such that
$F$ maps objects of $\frak J$ to objects of $\frak J'$.
Then $\frak J$ is $(F'\circ F)$-injective and we have
a natural isomorphism:
$$
  R(F'\circ F) \ = \ RF'\circ RF\ .
$$
\endproclaim

For the proof of this result see, for example, \cite{KSa, Prop. 1.8.7}

%
%

\head {Cohomology of sheaves} \endhead

\noindent
Let $X$ be a topological space. Let $\frak C$ be the category of $\CC$-vector spaces 
and let $\frak C(X)$ be the category  of sheaves of $\CC$-vector spaces on $X$.

\definition{Definition}
Let $\Gamma(X,-)$ be the global section functor from $\frak C(X)$
to $\frak C$. This functor is left exact. We denote by $\oh^i(X,-)$
the $i$-th derived functor $R^i\Gamma(X,-)\,$.
For a given sheaf $\F$ on $X$, the vector space
$$
  \oh^i(X,\F)\ = \ R^i@!@!\Gamma(X,\F)
$$
is called the $i$-{\it th cohomology space\/} 
of $X$ with coefficients in the sheaf $\F$.
\enddefinition

\definition{Injective, flabby, and c-soft resolutions}
Recall that to compute $R\Gamma(X,\F)$, we replace
$\F$ by a complex of injective sheaves $\Cal I^\cdot$ quasi-isomorphic
to $\F$ and then apply the functor $\Gamma(X,-)$
to $\Cal I^\cdot$. It is actually enough to choose the $\Cal I^i$
in a subcategory  which is injective with respect to the functor
$\Gamma(X,-)$.
Examples of such categories are the category of flabby sheaves, 
and in the case when $X$ is locally compact, the category of c-soft sheaves.
\enddefinition

\remark{Example}
Let $X$ be a real $\Cal C^\infty$-manifold of dimension $n$ and
let $\Cal E^{p}$ be the sheaf of smooth
$p$-forms on $X$. The sheaves $\Cal E^{p}$ are $c$-soft.
By the Poincar\'e Lemma, the {\it de Rham complex\/}
$$
  0 @>>> \CC_X @>>> \Cal E^{0} @>d>> \cdots @>>>\Cal E^{n} @>>> 0
$$
is exact. Thus the constant sheaf $\CC_X$ is quasi-isomorphic
to the complex $\Cal E^\cdot$, and for any $p$ the cohomology space 
$\oh^p(X,\CC_X)$ is the space of globally closed $p$-forms modulo
the space of globally exact $p$-forms.
\endremark

\remark{Axiomatic sheaf cohomology}
Let $\frak J(X)$ be any full additive subcategory of $\frak C(X)$ which
is injective with respect to $\Gamma(X,-)\,$.
Then the functors $\oh^i: \frak
C(X) @>>> \frak C\,$ satisfy the following properties: 
\vskip 5pt
\roster 
\item"(i)" There is a natural isomorphism $\Gamma(X,-) = \oh^0(X,-)\,$.
\item"(ii)" If $\, 0 @>>> \F' @>>> \F @>>> \F'' @>>> 0\,$ is a short exact
sequence of sheaves then there is a long exact sequence
$$
 \quad\quad \cdots @>>> \oh^{i}(X, \F') @>>> \oh^i(X,\F)
  @>>> \oh^i(X,\F'') @>{\partial^i}>> \oh^{i+1}(X, \F')  @>>> \cdots \ ,
$$
and the connecting homomorphisms $\partial^i$ behave functorially.
\item"(iii)" If $\Cal J$ is any sheaf of the category $\frak J(X)$, then
$\oh^i(X,\Cal J)=0\,$ for $i\not = 0\,$.
\endroster
\vskip 5pt\noindent
The functors $\oh^i$ are uniquely determined by these axioms. 
\endremark


\head Local cohomology\endhead

\definition{Sections with supports}
Let $Z$ be a locally closed subset of $X$.
We choose an open subset $V$ of $X$ containing $Z$ as a closed subset, and 
then define 
$$
  \Gamma_Z(X,\F)\ = \ \{ s\in \F(V)\ :\ s@,@,|_{V- Z}\,=\,0\,\}\ .
$$
One checks that $\Gamma_Z(X,\F)$ is independent of the choice of the open 
subset $V$. We call $\Gamma_Z(X,\F)$ the {\it sections of $\F$ with support
in $Z$\/}. If $Z=Y$ is closed then $\Gamma_Y(X,\F)$ is just the global sections
with support in $Y$. If $\,Z=V$ is open then $\Gamma_V(X,\F) = \F(V)$.
\vskip 5 pt
Let $X,Z$, and $V$ be as above. Then if $U$ is an open subset of $X$,
the natural restriction homomorphism $\F(V) @>>> \F(V\cap U)$
induces a homomorphism $\Gamma_Z(X,\F) @>>> \Gamma_{Z\cap U}(U,\F|_U)$.
The presheaf $U\mapsto \Gamma_{Z\cap U}(U,\F|_U)$ is a sheaf. This
sheaf is denoted by $\Gamma_Z(\F)$, and is called the {\it sheaf of sections
of $\F$ with support in $Z$}. 
\vskip 5 pt
The functors $\Gamma_Z(X,-): \frak C(X) @>>>  \frak C$ and $\Gamma_Z: \frak C(X)
@>>>  \frak C(X)$ are left exact. Moreover, we have
$\Gamma_Z(X,-) = \Gamma(X,-) \circ \Gamma_Z(-)$.
\enddefinition

\remark{Remark}
There is a different interpretation of the functors above as follows.
Let $\CC_Z$ be the constant sheaf on $Z$, which we also may interpret as a sheaf
on $X$ by extending it by zero outside $Z$.
Then we have natural isomorphisms of functors
$\Gamma_Z(X,-) = \Hom(\CC_Z,-)\,$, and $\Gamma_Z(-) = \ShHom(\CC_Z,-)$.
\endremark

\definition{Definition}
The $i$-th right derived functors of $\Gamma_Z(X,-)$ and $\Gamma_Z$ are 
by denoted by $\oh^i_Z(X,-)$ and $\oh^i_Z(-)$, respectively.
For a given sheaf $\F$, the vector space $\oh^i_Z(X,\F)$ (resp., the sheaf
$\oh^i_Z(\F)\,$) is called the {\it $i$-th cohomology space (resp., cohomology
sheaf) of $X$ with coefficients in $\F$ and supports in $Z$\/}. Note that
if $Z=X$ then  $\oh^i_Z(X,\F)=  \oh^i(X,\F)$ is the usual sheaf cohomology.
The natural properties of the functors $\oh^i(X,-)$ generalize to properties
of the functors $\oh^i_Z(X,-)$.
\enddefinition 

\remark{Remark}
For any sheaf $\F$ we have a canonical isomorphism
$$
  \oh^i_Z(X,\F)\ = \ \Ext^i(\CC_Z,\F) \ ,
$$
and similarly, $\oh^i_Z(\F) =  \ShExt^i(\CC_Z,\F)$, where $\ShExt$ 
is the derived functor obtained from $\ShHom$.
\endremark

\definition{Relative cohomology}
We often also write
$$
  \oh_Z^i(X,\F)\ = \ \oh^i(X,X-Z;\F)\ ,
$$
and think of the cohomology spaces $\oh_Y^i(X,\F)$ as {\it relative cohomology
of the pair $(X,X-Z)$ with coefficients in $\F$\/}. 
\enddefinition

\definition{Long exact sequences}
Let $Z$ and $X$ be as above, and let $Y$ be closed in $Z$.
Then for any sheaf $\F$, there is an exact sequence
$\,0 @>>> \Gamma_Y(\F) @>>> \Gamma_Z(\F) @>>> \Gamma_{Z-Y}(\F)$.
Moreover, if $\F$ is flabby, then this sequence extends to a short exact sequence. 
Hence we get long exact sequences
$$
  \cdots @>>> \oh^i_Y(\F) @>>> \oh^i_Z(\F)
@>>> \oh^{i}_{Z-Y}(\F) @>>> \oh^{i+1}_Y(\F) @>>>
\cdots
$$
\vskip -5 pt \noindent
and
\vskip -10 pt
$$
  \cdots @>>> \oh^i_Y(X,\F) @>>> \oh^i_Z(X,\F) @>>> 
  \oh^{i}_{Z-Y}(X,\F)  @>>> \oh^{i+1}_Y(X,\F) @>>> \cdots \ .
$$
\enddefinition

\definition{Excision}
Let $Z$ be a locally closed subset of $X$,  and let $V$ be an open subset of $X$
containing $Z$.  Then for any sheaf $\F$, there exists a natural isomorphism
$$
  \oh^i_Z(X,\F)\ = \ \oh^i_Z(V,\F|_V)\ .
$$
\enddefinition


\Bibliography			
\widestnumber\key{MMMM}		

\ref
\key{\bf BV} 
\by D.Barbasch and D.Vogan 
\paper The local structure of characters 
\jour Jour. Func. Anal. 
\vol 37 
\yr 1980 
\pages 27--55
\endref

\ref
\key{\bf BB1} 
\by A. Beilinson and J. Bernstein 
\paper Localisations de $\frak g$--modules 
\pages 15--18
\yr 1981 
\vol 292 
\jour C. R. Acad. Sci. Paris
\endref

\ref
\key{\bf BB2} 
\by A. Beilinson and J. Bernstein 
\paper A generalization of Casselman's submodule theorem
\pages  35--52
\inbook in: Representation Theory of Reductive Groups,
 Progress in Mathematics, Vol 40, Birk\-h\"au\-ser, Boston, 1983
\endref

\ref
\key {\bf BBD}
\by A. Beilinson, J. Bernstein, and P. Deligne
\paper Faiseaux pervers
\jour Ast\'erisque
\vol 100
\yr 1983
\endref

\ref
\key{\bf BL} 
\by J. Bernstein and V. Lunts
\book Equivariant sheaves and functors,  {\rm Lecture Notes in Mathematics
1578}
\publ Springer
\yr 1994
\endref

\ref
\key{\bf BM} 
\by E. Bierstone and P.D. Milman
\paper Semi-analytic and subanalytic sets
\jour Publications Math. IHES
\vol 67
\yr 1988
\pages 5--42
\endref

\ref
\key{\bf Bo} 
\by A. Borel et al.
\book Algebraic D-Modules
\bookinfo Perspectives in Mathematics
\publ Academic Press
\yr 1987
\endref

\ref
\key{\bf C} 
\by J.-T. Chang 
\paper Asymptotics and characteristic cycles for representations of complex groups 
\jour  Compositio Math.
\yr 1993
\pages 265 -- 283
\endref

\ref
\key {\bf DM}
\by  L. van den Dries and C. Miller
\paper Geometric categories and o-minimal structures
\jour Duke Math. Jour.
\vol 84
\yr 1996
\pages 497 -- 540
\endref

\ref
\key {\bf EW}
\by  T. J. Enright and N. R. Wallach
\paper Notes on homological algebra and
representations of Lie algebras
\jour Duke Math. Jour.
\vol 47
\yr 1980
\pages 1--15
\endref

\ref 
\key{\bf K1}
\by M. Kashiwara 
\paper Index theorem for constructible sheaves
\inbook Syst\`emes differentiels et singularit\'es
\eds A.Galligo, M.Maisonobe, and Ph. Granger
\jour Ast\'erisque \vol 130 \yr 1985 \pages 193 -- 209
\endref

\ref 
\key{\bf K2}
\by M. Kashiwara 
\paper Open problems in group representation theory
\inbook Proceedings of Tanig\-uchi Symposium held in 1986,  RIMS preprint 569
\publ Kyoto University
\yr 1987
\endref

\ref 
\key{\bf KSa} 
\by M. Kashiwara and P. Schapira 
\book Sheaves on manifolds
\publ Springer \yr 1990 \endref

\ref 
\key {\bf KSd} 
\by M.Kashiwara and W.Schmid 
\paper Quasi-equivariant $\Cal D$-modules, equivariant derived category,  and
representations of reductive Lie groups
\inbook Lie Theory and Geometry, in Honor of Bertram Kostant 
\bookinfo {\rm Progress in Mathematics}
\publ Birkh\"auser
\yr 1994
\pages 457--488
\endref

\ref
\key {\bf MP} 
\by D. Mili\v ci\'c and P. Pand\v zi\'c
\paper Equivariant derived categories,
Zuckerman functors and localization
\inbook  Geometry and representation theory of real and
$p$-adic groups (C\'ordoba, 1995)
\bookinfo {\rm  Progr. Math., 158}
\publ Birkh\"auser
\yr 1998 
\pages 209--242
\endref

\ref
\key{\bf MUV} 
\by I. Mirkovi\'c, T. Uzawa, and K. Vilonen 
\paper Matsuki correspondence for sheaves
\jour Inventiones Math. 
\vol 109 
\yr 1992 
\pages 231--245
\endref

\ref
\key {\bf MV} 
\by I. Mirkovi\'c and K. Vilonen 
\paper Characteristic varieties of character sheaves 
\jour Inventiones Math. 
\vol93 
\yr 1988 
\pages 405--418
\endref

\ref
\key{\bf N}
\by L. Ness
\paper A stratification of the null cone via the moment map
\jour Am. Jour. Math. \vol 106 \yr 1984 \pages 1281 -- 1325
\endref

\ref
\key{\bf R1} 
\by W. Rossmann 
\paper Kirillov's character formula for 
reductive Lie groups
\jour Inventiones Math. \vol 48 \yr 1978 \pages 207--220
\endref

\ref
\key {\bf R2} 
\by W. Rossmann 
\paper Characters as contour integrals
\inbook  Lecture Notes in Mathematics 1077
\publ Springer
\pages 375--388
\yr1984
\endref

\ref
\key{\bf S1} 
\by W. Schmid 
\paper Boundary value problems for group invariant differential equations
\inbook The mathematical heritage of \`Elie Cartan
\jour Ast\'e\-risque, Numero Hors Serie 
\yr 1985
\pages 311--321
\endref

\ref
\key{\bf S2} 
\by W. Schmid 
\paper Construction and classification of irreducible Harish-Chandra
modules
\inbook Harmonic analysis on reductive groups (Brunswick, ME, 1989) 
\bookinfo Progr. Math., 101 
\publ Birkh\"auser 
\yr1991.
\pages 235--275
\endref

\ref
\key{\bf SV1} 
\by  W. Schmid and K. Vilonen 
\paper Characteristic cycles of constructible sheaves
\jour Inventiones Math.
\vol 124
\yr 1996
\pages 451--502
\endref

\ref
\key{\bf SV2} 
\by  W. Schmid and K. Vilonen 
\paper Two geometric character formulas for reductive Lie groups
\jour Jour. AMS
\yr 1998
\vol 11
\pages 799 -- 876
\endref

\ref
\key{\bf SV3} 
\by  W. Schmid and K. Vilonen 
\paper Characteristic cycles and wave front cycles of representations of reductive
Lie groups
\paperinfo to appear in the Annals of Math
\endref

\ref
\key{\bf SV4} 
\by  W. Schmid and K. Vilonen 
\paper On the geometry of nilpotent orbits
\paperinfo to appear in the Atiyah volume of Asian Journal of Mathematics
\endref

\ref
\key{\bf SW} 
\by W. Schmid and J. Wolf 
\paper Geometric quantization and
derived functor modules for semi\-simple Lie groups 
\jour Jour. Funct. Anal.
\vol 90 
\yr 1990
\pages 48--112
\endref

\enddocument

%% file: macros.tex
\define\name#1{{\smc #1\/}}
\define\db#1{\operatorname {D}^b(#1)}
\define\od{\operatorname {D}}
\define\occ{\operatorname {CC}}
\define\cc{\operatorname {CC}}
\define\ct{T^*X}
\define\oh{{\operatorname H}}
\define\ok{{\operatorname K}}
\define\ohf#1{{\oh^{inf}_{#1}}}
\define\oc#1{{\operatorname C_{#1}^{inf}}}
\define\oz#1{{\operatorname Z_{#1}^{inf}}}

\define\fg{{\frak g}}

\define\fk{{\frak k}}
\define\fb{{\frak b}}

\define\fn{{\frak n}}
\define\fh{{\frak h}}

\define\ft{{\frak t}}
\define\fp{{\frak p}}

\define\bc{{\Bbb C}}

\define\TR{{T_\Bbb R}}

\define\GR{{G_\Bbb R}}
\define\UR{{U_\Bbb R}}
\define\KR{{K_\Bbb R}}
\define\gr{{\fg_{\Bbb R}}}
\define\kr{{\fk_{\Bbb R}}}

\define\ug{{\Cal U(\fg)}}
\define\zg{{\Cal Z(\fg)}}
\define\cf{{\Cal F}}
\define\cg{{\Cal G}}
\redefine\O{{\Cal O}}

\define\h{\hat}
\define\TP{{\Theta_\pi}}
\define\tp{{\theta_\pi}}

\define\Ext{\operatorname{Ext}}
\define\Hom{\operatorname{Hom}}

\redefine\Re{\operatorname{Re}}

\define\CC{{\Bbb C}}
\define\RR{{\Bbb R}}

\predefine\greater{\gg}
\redefine\gg{{\frak g}}

\define\ns#1{\Vert {#1}\Vert^2}


\redefine\ker{\operatorname{ker}}

\define\im{\operatorname{im}}

\define\Ad{\operatorname{Ad}}

\define\ZZ{{\Bbb Z}}
\define\G{{\Cal G}}
\define\F{{\Cal F}}
\define\T{{\Cal T}} 

\define\id{\operatorname{id}}
\define\ShExt{\operatorname{{\Cal E}\!@,@,@,\text{\it x\/}@!\text{\it
t\/}}}
\define\ShHom{\operatorname{{\Cal H}\!\text{\it o}@!@!@!\text{\it
m}\,@!@!@!}}
\define\RShHom{R\!@,@,\operatorname{{\Cal H}\!\text{\it o}@!@!@!\text{\it
m}\,@!@!@!}}
\define\RHom{R\!@,\operatorname{Hom}}
\define\RGamma{R@,@,\Gamma}
\define\Dual{\Bbb D}
\define\pt{\text{\rm pt}}
\redefine\phi{\varphi}
\define\supp{\operatorname{supp}}
\define\Ind{\operatorname{Ind}}
\define\qis{@>{}_{qis}>>}
\define\leftqis{@<{}_{qis}<<}